\documentclass[numbook, envcountsame, smallcondensed]{svjour3} 

\usepackage{cite} 
\usepackage{amsmath,amscd}
\usepackage{amssymb} 
\usepackage[retainorgcmds]{IEEEtrantools}
\usepackage{graphicx}
 \usepackage{tikz-cd}
\usepackage{comment}
\usepackage[polish, english]{babel}
\usepackage{enumerate}
\usepackage{mathtools}
\usepackage{yfonts}

\def\makeheadbox{\relax}
\DeclarePairedDelimiter{\norm}{\lVert}{\rVert}

\def\rd{\mathrm{d}}
\DeclarePairedDelimiter\abs{\lvert}{\rvert}

\begin{document}
\title{Differentiable functions on modules and equation $grad(w)=\tens{M}grad(v)$}

\author{Krzysztof J. Ciosmak}

\institute{Krzysztof J. Ciosmak \at University of Oxford, Mathematical Institute,\\
Andrew Wiles Building, Radcliffe Observatory Quarter,\\
Woodstock Rd, Oxford OX2 6GG, United Kingdom,\\
\email{ciosmak@maths.ox.ac.uk}, \and 
University of Oxford, St John's College,\\
St Giles', Oxford OX1 3JP,United Kingdom,\\
\email{krzysztof.ciosmak@sjc.ox.ac.uk}.
}
\label{firstpage}
\maketitle
\begin{abstract}

Let $A$ be a finite-dimensional, commutative algebra over $\mathbb{R}$ or $\mathbb{C}$. We extend the notion of $A$-differentiable functions on $A$ and develop a theory of $A$-differentiable functions on finitely generated $A$-modules. Let $U$ be an open, bounded and convex subset of such a module. We give an explicit formula for $A$-differentiable functions on $U$ of prescribed class of differentiability in terms of real or complex differentiable functions, in the case when $A$ is singly generated and the module is arbitrary and in the case when $A$ is arbitrary and the module is free. We prove that certain components of $A$-differentiable function are of higher differentiability than the function itself. 

Let $\tens{M}$ be a constant, square matrix. Using the aforementioned formula we find a complete description of solutions of the equation $grad(w)=\tens{M}grad(v)$.

We formulate the boundary value problem for generalised Laplace equations $\tens{M}\nabla^2 v=\nabla^2v \tens{M}^{\intercal}$ and prove that for given boundary data there exists a unique solution, for which we provide a formula.
\end{abstract}

\begin{acknowledgements}
I would like to express my gratitude towards Micha\l{} Wojciechowski, advisor of my Master's thesis, from which this work has emerged. \\
This is a preprint of the work accepted for publication in the journal \emph{Algebra and Analysis}, 2022. The owner of distribution rights - POMI RAS.
\end{acknowledgements}

\section{Introduction}

Let $\tens{M}$ be a square $n\times n$ matrix. Consider a system of partial differential equations
\begin{equation}\label{eqn:e}
grad(w)=\tens{M}grad(v),
\end{equation}
where $v, w$ are real-valued functions defined on an open set $U\subset \mathbb{R}^n$. The first systematic treatment of this system and a complete description of solutions on convex sets was provided in \cite{JO} by Jodeit and Olver. Later, Waterhouse in \cite{W1} observed that if $\tens{M}$ has only one Jordan block corresponding to each eigenvalue, then any solution to the equation is represented by some differentiable function on a finite-dimensional, commutative algebra $A$. Differentiable functions on algebras (see \cite{K, R, W2}), studied since their definition in 1891 by Georg Scheffers -- a student of Sophus Lie (see \cite{W2} for more historical background), are not enough to provide such representation for an arbitrary $\tens{M}$. However, Waterhouse in \cite{W1} suggested that a solution in a general case could be represented by an $A$-differentiable function on a module. Following this suggestion we develop the theory of $A$-differentiable functions on modules (see \S \ref{S:Def}). We prove that for an arbitrary $\tens{M}$ any solution may be represented as an $A$-differentiable function on a certain $A$-module (see Theorem \ref{T:pairs}).  

Waterhouse in \cite{W2} showed that $A$-differentiable functions on algebras depend polynomially on some coordinates. We extend this result and provide (see \S\ref{S:structure}) an exact dependence of such functions on these coordinates and give a formula, in terms of real and complex differentiable functions, for any $A$-differentiable function. In \S \ref{S:structure} we investigate the structure of any differentiable function on a free module (see Theorem \ref{T:structurefree}) and the structure of any differentiable function on an arbitrary module over a singly generated algebra (see Theorem \ref{T:structure-one}). Surprisingly, even in the case of $\mathbb{R}$-algebras, certain components of an $A$-differentiable function are of the higher differentiability than it is assumed of the function. 

Jonasson in \cite{J1} studied the equation (\ref{eqn:e}) and posed a question, whether it is true that every solution is given by a power series. We show that the answer depends on the roots of characteristic polynomial of $\tens{M}$. If there are no real roots, then the answer is positive and, conversely, if the opposite case occurs, then the answer is negative (see Theorem \ref{T:analyticity} and \S \ref{S:analiticity}).

When analysing the equation (\ref{eqn:e}), the generalised Laplace equations 
\begin{equation}\label{eqn:gle}
\tens{M}\nabla^2 v=\nabla^2v \tens{M}^{\intercal}
\end{equation}
appear naturally, as these are exactly the integrability conditions for (\ref{eqn:e}). We are thus interested in the boundary value problem for equations (\ref{eqn:gle}). \S \ref{S:boundary} is devoted to this problem. The connection with $A$-differentiable functions on modules (see Theorem \ref{T:components}) and the structure theorem (Theorem \ref{T:structurefree} and Theorem \ref{T:structure-one}) allows us to provide an appropriate formulation of boundary value problem (see Theorem \ref{T:boundary}) such that (\ref{eqn:gle}) has a unique solution for a given boundary data. Moreover, there exists a continuous linear operator that maps boundary data to the unique solution.

There are several questions which arise. Let us enumerate some of them.
\begin{enumerate}[\upshape (i)]
\item How does an $A$-differentiable function between arbitrary $A$-modules look?  
\item How to generalise the notion of $A$-differentiability to non-commutative algebras so that any $A$-analytic function is $A$-differentiable and the set of $A$-differentiable functions forms an algebra?
\item How to solve the equation $grad(w)=\tens{M}grad(v)$ when $\tens{M}$ is a non-constant matrix?
\item What are other partial differential equations such that the space of solutions forms an algebra? What are the conditions for such a property? When is such an algebra semisimple? When is it finitely generated?
\end{enumerate}

\section{Differentiable functions on modules}\label{S:Def}

We begin with a definition of differentiability of functions on modules. Throughout the paper $A$ denotes a finite-dimensional \emph{commutative} $\mathbb{F}$-algebra with unit over field $\mathbb{F}\in \{\mathbb{R},\mathbb{C}\}$. $B, C$ stand for finitely generated $A$-modules. $U\subset B$ is an open set in the norm topology of $B$. We shall denote by $\norm{\cdot}$ a norm. Usually it can be chosen arbitrarily, otherwise we shall specify what conditions we want a norm to satisfy.

Recall that a module $B$ over $A$ is finitely generated if there exists a surjective module homomorphism $\eta \colon A^n \to B$. For finite-dimensional $A$ this happens if and only if $B$ is finite-dimensional as a vector space. 

If $B$ is a finitely generated free module, we may choose an $A$-basis $b_1,\dotsc,b_n$ for $B$. For 
\begin{equation*}
x=\sum_{i=1}^nx_ib_i\text{, }y=\sum_{i=1}^ny_ib_i\text{ and a matrix }\tens{C}=[c_{ij}]_{i,j=1,\dotsc ,n}\in M_{n\times n}(A)
\end{equation*}
we write 
\begin{equation*}
\langle x,y\rangle_A=\sum_{i=1}^{n}x_iy_i\text{ and }\tens{C}y=\sum_{j=1}^{n}c_{ij}y_{j}b_i.
\end{equation*}

\begin{definition}\label{D:A-differentiability}
A function $f\colon U \to C$ is said to be $A$-\emph{differentiable} if there is a map $\mathrm{D}f\colon U\times B \to C$, such that
\begin{enumerate}[\upshape (i)]
\item $\lim_{y\to0}\tfrac{1}{\norm{y}}\norm{f(x+y)-f(x)-\mathrm{D}f(x)(y)}=0$,
\item $\mathrm{D}f(x)(ay)=a\mathrm{D}f(x)(y)$ for all $x\in U$, $y\in B$ and $a\in A$.
\end{enumerate}

\end{definition}

\begin{remark}
If $f$ is differentiable as a function of its real variables and the linear derivative $\mathrm{D}f\colon U\times B \to C$ satisfies $\mathrm{D}f(x)(ay)=a\mathrm{D}f(x)(y)$ for all $x\in U$, $y\in B$ and $a\in A$, then $f$ is $A$-differentiable. 
\end{remark}

\begin{remark}
As $B, C$ and $A$ are finite-dimensional, all norms on these spaces are equivalent. Thus the condition of $A$-differentiability does not depend on the choice of the norm.
\end{remark}

Let us recall the definition of differentiability on algebras  (see \cite{W1}, \cite{W2}) -- a function $f\colon U\to A$ on $U\subset A$ is said to be $A$-differentiable if for every $x\in U$ there exists a limit
\begin{equation*}
\lim_{{\substack{y\to x\\
                  y\in (x+ G(A))\cap U}}} \frac{f(y)-f(x)}{y-x},
\end{equation*}
where $G(A)$ denotes the group of invertible elements in $A$.
The above definition of $A$-differentiability agrees with Definition \ref{D:A-differentiability} if we put $B=C=A$. Indeed, as is proved in \cite{W2}, if the above limit exists, then $f$ is differentiable as a function of its real variables. The derivative at a point $x\in U$ is given by a multiplication by a certain element $f'(x)\in A$, so it is $A$-linear. The converse is trivial. 

We shall indicate that the condition stated in the definition constitutes an equivalent set of partial differential equations for $f$ seen as a function of real variables. We write down these equations for $A$-valued functions. Let $c_1,\dotsc,c_t$ denote an $\mathbb{F}$-basis of $B$ for the moment. Let $e_1,\dotsc,e_k$ denote an $\mathbb{F}$-basis of $A$, and let 
\begin{equation*}
e_ic_j=\sum_{s=1}^{t} \alpha_{ij}^{s}c_s\text{ for some }\alpha_{ij}^s\in\mathbb{F}.
\end{equation*}
The condition for $A$-differentiability is equivalent to the equations
\begin{equation*}
\sum_{s=1}^{t}\alpha_{ij}^{s}\mathrm{D}f(\cdot)(c_s)=e_i\mathrm{D}f(\cdot)(c_j)\text{ for all }i=1,\dotsc,k\text{ and all }j=1,\dotsc t.
\end{equation*}
These equations are reduced to the Cauchy--Riemann equations when $A=\mathbb{C}$, viewed as an $\mathbb{R}$-algebra, and $B=A$. When $B=\mathbb{C}^n$ and $A=\mathbb{C}$, then these are the conditions for a function of several complex variables to be holomorphic.

\begin{lemma}\label{L:secondderivative}
Let $k\geq 1$. Let $f\colon U\to C$ be a $C^k$ function that is $A$-differentiable. Then for any $x_1,\dotsc,x_{k-1}\in B$ the function $\mathrm{D}^{k-1}f(\cdot)(x_1,\dotsc,x_{k-1})\colon U\to C$ is again $A$-differentiable. Moreover 
\begin{equation*}
\mathrm{D}^kf(\cdot)(ax_1,\dotsc,x_k)=a\mathrm{D}^kf(\cdot)(x_1,\dotsc,x_k)
\end{equation*}
for all $x_1,\dotsc,x_k\in B$ and all $a\in A$.
\end{lemma}
\begin{proof} 
Let $g(\cdot)=\mathrm{D}^{k-1}f(\cdot)(x_1,\dotsc,x_{k-1})$. We have 
\begin{equation*}
\begin{aligned}
\mathrm{D}^kf(\cdot)(ax_1,\dotsc,x_k)=&\mathrm{D}^{k-1}(\mathrm{D}f(\cdot)(ax_1))(x_2,\dotsc,x_k)=\\
=&\mathrm{D}^{k-1}(a\mathrm{D}f(\cdot)(x_1))(x_2,\dotsc,x_k)=a\mathrm{D}^kf(\cdot)(x_1,\dotsc,x_k).
\end{aligned}
\end{equation*}
Thus $\mathrm{D}g(\cdot)(ay)=a\mathrm{D}g(\cdot)(y)$.
\end{proof}

\begin{lemma}\label{L:etadifferentiability}
Let $B, C, D$ be $A$-modules and let $U\subset B$ be an open set. Let $f\colon U \to D$. Let $\pi\colon C\to B$ be an $A$-linear surjection. Then
\begin{enumerate}[\upshape (i)]
\item $f$ is $A$-differentiable on $U$ if and only if $f\circ \pi$ is $A$-differentiable on $\pi^{-1}(U)$,
\item $f$ is $k$-times continuously differentiable on $U$ if and only if $f\circ \pi $ is $k$-times continuously differentiable on $\pi^{-1}(U)$,
\item $f$ is $k$-times continuously differentiable on $\overline{U}$ if and only if $f\circ \pi $ is $k$-times continuously differentiable on $\overline{\pi^{-1}(U)}$,
\end{enumerate}
\end{lemma}
\begin{proof}
One needs to verify that $\mathrm{D}f$ exists and is $A$-linear if $\mathrm{D}(f\circ \pi)$ exists and is $A$-linear. Let us first show that $\mathrm{D}f$ exists. 
 
For some linear subspace $V\subset C$, we have $C=V \oplus\ker \pi $, as vector spaces. On $C$ we may introduce a norm by the formula
\begin{equation*}
\norm{r+s}=\norm{\pi(r)}_B+\norm{s}_{C}
\end{equation*}
for $r\in V$ and $s\in \ker\pi$ and some norms $\norm{\cdot}_B$ on $B$ and $\norm{\cdot}_C$ on $\ker \pi$.

For any given $x\in U$, $y\in B$ choose $t\in\pi^{-1}(x)$ and $r\in \pi^{-1}(y)\cap V$ such that $\norm{y}_B= \norm{\pi(r)}_B=\norm{r}$. Since $\norm{r}\to 0$ if $\norm{y}_B\to 0$, we have
\begin{equation*}
\begin{aligned}
&\tfrac{1}{\norm{y}_B}\norm{f(x+y)-f(x)-\mathrm{D}(f\circ \pi)(t)(r)}=\\
&=\tfrac{1}{\norm{r}}\norm{f(\pi(t+r))-f(\pi(t))-\mathrm{D}(f\circ \pi)(t)(r)}\xrightarrow[{y\to 0}]{}0.
\end{aligned}
\end{equation*}
Hence $f$ is differentiable and $\mathrm{D}f(x)(y)=\mathrm{D}(f\circ \pi)(t)(r)$.
$A$-linearity follows from $A$-linearity of $\mathrm{D}(f\circ \pi)$ and $\pi$. Indeed, for $a\in A$ we write $ar=v_a+k_a$, with $v_a\in V$, $k_a\in \ker\pi$. Then $\pi(ar)=\pi(v_a)$ so $v_a\in \pi^{-1}(ay)\cap V$. Thus
\begin{equation*}
\begin{aligned}
&\mathrm{D}f(x)(ay)=\mathrm{D}f(\pi(t))(\pi(v_a))=\mathrm{D}(f\circ \pi)(t)(v_a)=\\
&=\mathrm{D}(f\circ \pi)(t)(ar-k_a)=a\mathrm{D}f(x)(y).
\end{aligned}
\end{equation*}

Let us prove that $f$ is $k$-times continuously differentiable if $f \circ \pi $ is $k$-times continuously differentiable. We shall proceed inductively. Assume that we have shown that for some $l\leq k-1$ we have
\begin{equation}\label{eqn:pidiff}
\mathrm{D}^lf(x)(y_1,\dotsc,y_l)=\mathrm{D}^l(f\circ \pi)(t)(r_1,\dotsc,r_l)
\end{equation}
for $t\in\pi^{-1}(x)\cap V$ and $r_i\in \pi^{-1}(y_i)\cap V$, $i=1,\dotsc,l$.
Then, by Lemma \ref{L:secondderivative}, for given $y_1,\dotsc, y_l$ function $\mathrm{D}^lf(\pi(\cdot))(y_1,\dotsc, y_l)$ is $A$-differentiable. The previous argument shows that $\mathrm{D}^lf(\cdot)(y_1,\dotsc, y_l)$ is $A$-differentiable and
\begin{equation*}
\begin{aligned}
\mathrm{D}^{l+1}f(x)(y_1,\dotsc,y_l, y_{l+1})=&\mathrm{D}(\mathrm{D}^l(f( \pi(\cdot))(y_1,\dotsc,y_l))(t)(r_{l+1})=\\
=&\mathrm{D}^{l+1}(f\circ \pi)(t)(r_1,\dotsc,r_{l+1})
\end{aligned}
\end{equation*}
for some $r_{l+1}\in \pi^{-1}(y_{l+1})\cap V$.
It follows that $f$ is $k$-times differentiable if $f\circ \pi$ is. 

Assume that $f\circ \pi$ is $k$-times continuously differentiable on $\pi^{-1}(U)$. Then, by (\ref{eqn:pidiff}), $f$ is $k$-times continuously differentiable on $u$. For if $(x_n)_{n=1}^{\infty}\subset U$ converges to $x$, then choosing 
\begin{equation*}
t_n\in \pi^{-1}(x_n)\cap V\text{ and } t\in \pi^{-1}(x)\cap V
\end{equation*}
we have $\norm{t_n-t}_B=\norm{\pi(x_n-x)}_B=\norm{x_n-x}$, so $t_n$ converges to $t$. 

Assume that $f\circ \pi $ is $k$-times continuously differentiable on $\overline{\pi^{-1}(U)}$. Then (\ref{eqn:pidiff}) allows us to extend the derivatives of $f$ continuously on $\overline{U}$, as we have $\pi^{-1}(\overline{U})=\overline{\pi^{-1}(U)}$.
\end{proof}

\subsection{Equivalent characterisations of $A$-differentiability}

We shall deal now with $A$-valued forms on $U$. As previously, we shall denote by $e_1,\dotsc ,e_k$ a basis of $A$ over $\mathbb{F}$. By $x_1,\dotsc,x_k$ we now denote the corresponding coordinate functions. If $x\colon U \to A$ is a joint $A$-coordinate function, then 
\begin{equation*}
x=\sum_{i=1}^k x_ie_i\text{ and }
\mathrm{d} x=\sum_{i=1}^k e_i \mathrm{d} x_i.
\end{equation*}

It is clear that the Poincar\'e's lemma and standard rules for differentiating wedge-product of two forms hold true.  

\begin{theorem}\label{T:equivalence-one}
Assume that $U\subset A$ is open and simply connected. Let $f\colon U\to A$ be continuously differentiable. The following conditions are equivalent
\begin{enumerate}[\upshape (i)]
\item\label{i:A-differentiability-one}$f$ is $A$-differentiable,
\item\label{i:closedness-one} the form $f \rd x$ is closed,
\item\label{i:integrals-one}  $\int_{\gamma}f \rd x=0$ for any smooth closed curve $\gamma$ in $U$,
\item\label{i:second derivative-one}there exists an $A$-differentiable function $g\colon U\to A$ that is twice continuously differentiable and such that 
\begin{equation*}
\mathrm{D}g(b)(x)=xf(b)
\end{equation*}
for all $b\in U$ and $x\in A$.
\end{enumerate}
\end{theorem}
\begin{proof}
Let us first see what the condition (\ref{i:closedness-one}) means. The form $f\rd x$ is closed if and only if $\rd x\wedge \rd f =0$. That is
\begin{equation*}
\begin{aligned}
0&=\left( \sum_{l=1}^{k}e_l \rd x_l\right)\wedge \left(\sum_{i=1}^{k}\mathrm{D}f(\cdot)(e_i)\rd x_i\right) =\\
&=\sum_{0\leq i<l\leq k}(e_l\mathrm{D}f(\cdot)(e_i)  - e_i \mathrm{D}f(\cdot)(e_l) )  \rd x_l \wedge \rd x_i.
\end{aligned}
\end{equation*}
This is equivalent to that for all $i,l=1,\dotsc,k$ we have $e_l \mathrm{D}f(\cdot)(e_i)=e_i \mathrm{D}f(\cdot)(e_l)$. By bilinearity of both sides we have equivalently 
\begin{equation}\label{eqn:A-commutativity-one}
y\mathrm{D}f(\cdot)(x)=x\mathrm{D}f(\cdot)(y)
\end{equation}
for all $x, y\in A$. Let $a\in A$. Then
\begin{equation}\label{eqn:A-linearity-one}
y\mathrm{D}f(\cdot)(ax)=ax\mathrm{D}f(\cdot)(y)=a(y\mathrm{D}f(\cdot)(x))=y(a\mathrm{D}f(\cdot)(x)).
\end{equation}
Letting $y$ to be the unit of $A$, we see that if (\ref{i:closedness-one}) is satisfied, then so is (\ref{i:A-differentiability-one}).
Conversely, knowing that (\ref{i:A-differentiability-one}) is true, then (\ref{eqn:A-linearity-one}) holds true and so does (\ref{eqn:A-commutativity-one}). This condition is equivalent to (\ref{i:closedness-one}).

If we know that there is an $A$-differentiable $g$ is as in (\ref{i:second derivative-one}), then by Lemma \ref{L:secondderivative}, we have
\begin{equation*}
x\mathrm{D}f(b)(ay)=\mathrm{D}^2g(b)(x,ay)=a\mathrm{D}^2g(b)(x,y)=x(a\mathrm{D}f(b)(y)),
\end{equation*}
thus $\mathrm{D}f(b)$ is $A$-linear.

Assume that the form $f\rd x$ is closed. By the Poincar\'{e}'s lemma there exists a twice continuously differentiable function $g\colon U \to A$ such that $\rd g=f\rd x$. That is
\begin{equation*}
\sum_{i=1}^k \mathrm{D}g(\cdot)(e_i)\rd x_i=\sum_{i=1}^k fe_i \rd x_i.
\end{equation*}
This condition is equivalent to $\mathrm{D}g(\cdot)(x)=xf(\cdot)$. Hence, $g$ is $A$-differentiable. 

Assume that condition (\ref{i:integrals}) is satisfied by $f\colon U\to A$. Choose any disc $D\subset U$, let $\gamma$ be its boundary. Then by the Stokes' theorem
\begin{equation*}
\int_D \rd f\wedge \rd x = \int_D \rd (f \rd x) =\int _{\gamma} f\rd x=0.
\end{equation*}
Since $D$ was arbitrary and $\rd f\wedge \rd x $ is continuous, we infer that (\ref{i:closedness}) holds true.

Assume that (\ref{i:closedness}) is fulfilled. Take any smooth closed curve $\gamma \colon [0,1]\to U$. By the Whitney's approximation theorem (see \cite{L}, Theorem 6.21) and by simple connectedness of $U$, $\gamma$ is a boundary of some smooth surface $D\subset U$. Again by the Stokes' theorem, the integral 
\begin{equation*}
\int_{\gamma} f \rd x=\int_{D} \rd(f \rd x )
\end{equation*}
vanishes. 
\end{proof}

Below, $B$ denotes a free module over $A$. If $U\subset B$ and $f\colon U \to B$ then we write $f=\sum_{i=1}^n f_i b_i$, where $b_1,\dotsc,b_n$ is an $A$-basis of $B$. We shall consider $A$-valued forms on $U$. In particular, let $z_i$ denote the $A$-valued $i^{\text{th}}$ coordinate function on $U$ and $\rd z_i$ the corresponding form. We write for $i=1,\dotsc,n$ further
\begin{equation*}
z_i=\sum_{r=1}^k x_{ir}e_r
\end{equation*}
with $x_{ir}$ for $r=1,\dotsc,k$ being the $\mathbb{F}$-coordinate functions.

\begin{proposition}\label{T:secondder}
Let $U\subset B$ be open and simply connected. Let $f\colon U\to B$ be continuously differentiable. The following are equivalent:
\begin{enumerate}[\upshape (i)]
\item $f$ is $A$-differentiable and such that $\mathrm{D}f(\cdot)$ satisfies 
\begin{equation*}
\mathrm{D}f_i(\cdot)(b_j)=\mathrm{D}f_j(\cdot)(b_i),
\end{equation*}
for all $i,j=1,\dotsc,n$,
\item\label{i:second} there exists a twice continuously $A$-differentiable function $g\colon U \to A$, such that 
\begin{equation*}
\mathrm{D}g(\cdot)(x)=\langle x,f(\cdot)\rangle_A.
\end{equation*}
\end{enumerate}
\end{proposition}
\begin{proof}
Define a one-form
\begin{equation*}
h=\sum_{i=1}^nf_i \rd z_i.
\end{equation*}
Then $\rd h=0$. Indeed, by $A$-linearity and the assumption we have
\begin{equation*}
\begin{aligned}
&\rd h=\sum_{i,s=1}^n\sum_{j,r=1}^k e_j\mathrm{D}f_i(\cdot)(e_rb_s) \rd x_{rs} \wedge \rd x_{ji}=
\sum_{i,s=1}^n \mathrm{D}f_i(\cdot)(b_s)\rd z_s\wedge \rd z_i=\\
&=\sum_{0\leq i<s\leq n}( \mathrm{D}f_i(\cdot)(b_s)-\mathrm{D}f_s(\cdot)(b_i)) \rd z_s\wedge\rd z_i=0.
\end{aligned}
\end{equation*}
By the Poincar\'{e}'s lemma there exists a twice continuously differentiable function $g\colon U \to A$ such that $\rd g=h$. That is 
\begin{equation*}
\sum_{i=1}^n\sum_{j=1}^k e_jf_i(\cdot)  \rd x_{ji}=\sum_{i=1}^n\sum_{j=1}^k \mathrm{D}g(\cdot)(e_jb_i) \rd x_{ji},
\end{equation*}
thus $\mathrm{D}g(\cdot)(e_jb_i)=e_jf_i(\cdot)$. Since this holds for all $j=1,\dotsc,k$, we have 
\begin{equation*}
\mathrm{D}g(\cdot)(yb_i)=yf_i(\cdot)
\end{equation*}
for all $y\in A$ and in consequence $\mathrm{D}g(\cdot)(x)=\langle x,f(\cdot)\rangle_A$ for all $x\in B$.

To prove the converse, let $g$ be as in (\ref{i:second}). Then $\mathrm{D}^2g(\cdot)(x,y)=\langle x, \mathrm{D}f(\cdot)(y)\rangle_A$. By Lemma \ref{L:secondderivative}, $\mathrm{D}f(\cdot)$ is $A$-linear. Further, since $\mathrm{D}^2g(\cdot)$ is symmetric, we have
\begin{equation*}
\langle x ,\mathrm{D}f(\cdot)(y)\rangle_A=\langle y,\mathrm{D}f(\cdot)(x)\rangle_A
\end{equation*}
Taking $x=b_i$ and $y=b_j$ we obtain $\mathrm{D}f_i(\cdot)(b_j)=\mathrm{D}f_j(\cdot)(b_i)$ for all numbers $i,j=1,\dotsc,n$.
\end{proof}

For any $b\in B$ define $t^b_j\colon A\to B$ by the formula
\begin{equation*}
t^b_j(a)=ab_j+b.
\end{equation*}
The following proposition tells, roughly speaking, that a differentiable function is $A$-differentiable if and only if it is $A$-differentiable on $A$-lines.

\begin{proposition}\label{P:separately}
Let $U\subset B$ be open. Let $f\colon U\to A$ be a differentiable function. The following conditions are equivalent:
\begin{enumerate}[\upshape (i)]
\item $f$ is $A$-differentiable,
\item for any $x\in B$ and $j=1,\dotsc,n$, 
\begin{equation*}
f\circ t^x_j\colon (t^x_j)^{-1}(U)\to A
\end{equation*}
is $A$-differentiable.
\end{enumerate}
\end{proposition}
\begin{proof}
If $f$ is $A$-differentiable, then $f\circ t^x_j$ is also $A$-differentiable, since $t^x_j$ is $A$-differentiable.   

Assume that $f\circ t^b_j$ is $A$-differentiable for any $b\in B$ and any $j=1,\dotsc,n$.
Choose $b\in U$ and $x\in B$, $x=\sum_{i=1}^na_ib_i$, and $a\in A$. Then
\begin{equation*}
\begin{aligned}
&\mathrm{D}f(b)(a\sum_{i=1}^n a_ib_i)=\sum_{i=1}^n \mathrm{D}f(b)(aa_ib_i)=\sum_{i=1}^n \mathrm{D}f(t^b_i(0))(\mathrm{D}t^b_i(0)(aa_i))\\
&=\sum_{i=1}^n \mathrm{D}(f\circ t^b_i)(0)(aa_i)=a\sum_{i=1}^n \mathrm{D}(f\circ t^b_i)(0)(a_i)=a\mathrm{D}f(b)(\sum_{i=1}^n a_ib_i).
\end{aligned}
\end{equation*}
This is to say, $Df(b)$ is $A$-linear.
\end{proof}

\begin{corollary}\label{C:equivalence}
Assume that $U\subset B$ is an open set such that for $j=1,\dotsc,n$ and any $b\in B$ the sets $(t^b_j)^{-1}(U)$ are simply connected.
Let $f\colon U \to B$ be a continuously differentiable function. The following conditions are equivalent
\begin{enumerate}[\upshape (i)]
\item\label{i:A-differentiability}$f$ is $A$-differentiable,
\item\label{i:closedness} for $j=1,\dotsc,n$ and any $b\in B$ the form $(f\circ t^b_j) \rd x$ is closed,
\item\label{i:integrals}  for $j=1,\dotsc,n$ and any $b\in B$ $\int_{\gamma}(f\circ t^b_j) \rd x=0$ for any smooth closed curve $\gamma$.
\end{enumerate}
\end{corollary}
\begin{proof}
Conditions (\ref{i:A-differentiability})-(\ref{i:integrals}) are equivalent by Theorem \ref{T:equivalence-one} and Proposition \ref{P:separately}.
\end{proof}

\begin{remark}
$A$-differentiability is a local condition. It it thus equivalent to the fact that conditions (\ref{i:closedness}) and (\ref{i:integrals}) hold locally. 
\end{remark}

\subsection{$A$-continuity}

The condition (\ref{i:integrals}) of Corollary \ref{C:equivalence} can be stated for merely continuous functions. For example, when $A=\mathbb{R}$, then it is satisfied by all continuous functions. If $B=A$, then it implies the existence of a primitive function, which is $A$-differentiable. Moreover, when $A$ is a $\mathbb{C}$-algebra, then every function which satisfies this condition is actually analytic. 

\begin{proposition}\label{P:A-continuous}
Assume that $U\subset A$ is an open, simply connected set. Let $f\colon U \to B$ be a continuous function. The following conditions are equivalent
\begin{enumerate}[\upshape (i)]
\item\label{i:A-continuity} $\int_{\gamma}f\rd x=0$ for any smooth closed curve $\gamma$,
\item\label{i:primitive}there exists a continuously $A$-differentiable function $g\colon U \to A$, such that $\mathrm{D}g(b)(x)=xf(b)$ for all $b\in U$ and $x\in B$.
\end{enumerate}
\end{proposition}
\begin{proof}
Choose a point $b_0\in U$. Let $b\in U$ and let $\gamma_b$ be a smooth curve with $\gamma_b(0)=b_0$ and $\gamma_b(1)=b$. Define $g(b)=\int_{\gamma_b}f\rd x$. Definition of $g$ is independent of the choice of smooth curve $\gamma_g$ as a consequence of the condition (\ref{i:A-continuity}). Then by the triangle inequality and by continuity of $f$ we have
\begin{equation*}
\begin{aligned}
&\lim_{x\to 0} \frac1{\norm{x}} \norm{g(b+x)-g(b)-xf(b)}=\\
&=\lim_{x\to 0} \frac1{\norm{x}} \Big\lVert \int_{[b,b+x]}f\rd x-xf(b)\Big\rVert\leq\\
&\leq \lim_{x\to 0} \frac1{\norm{x}} \int_0^1\norm{(f(b+tx)-f(b))x}\rd t=0.
\end{aligned}
\end{equation*}
This means exactly that $\mathrm{D}g(b)(x)=xf(b)$. 

Conversely, if (\ref{i:primitive}) holds then, as in Theorem \ref{T:equivalence-one}, $g$ satisfies $\rd g=f\rd x$. Choose a smooth curve $\gamma$. Again as in Theorem \ref{T:equivalence-one} we may assume that $\gamma$ is a boundary of a smooth surface $D$ and thus
\begin{equation*}
\int_{\gamma}f\rd x=\int_{\gamma}\rd g=0.
\end{equation*}
\end{proof}

\begin{definition}\label{D:A-continuity}
Let $A^n$ be a free module and $B$ be a finitely generated $A$-module, let $\eta\colon A^n\to B$ be a surjective $A$-linear map. Let $U\subset B$ be an open set. A function $f\colon U \to A$ is called $A$-\emph{continuous}, if it is continuous and such that every point $b$ in $U$ admits a neighbourhood $V\subset U$ such that 
\begin{equation*}
\int_{\gamma}(f\circ \eta \circ t^b_j)\rd x =0
\end{equation*}
for any smooth closed curve $\gamma\subset (\eta\circ t^b_j)^{-1}(V)$, any $j=1,\dotsc,n$.
\end{definition}

\begin{remark}
Assume that $f$ is continuously differentiable. Then Lemma \ref{L:etadifferentiability} and Corollary \ref{C:equivalence} imply that $f\colon U \to A$ is $A$-continuous if and only if it is $A$-differentiable. 
\end{remark}

\section{Components of $A$-differentiable functions}

We will here characterise $\mathbb{F}$-valued functions defined on an arbitrary finitely generated $A$-module $B$ that arise as components of $A$-differentiable functions on $B$. That is, we shall specify differential conditions for a function $v\colon U\to \mathbb{F}$ which are satisfied if and only if $v=\psi(f)$ for some $\mathbb{F}$-linear functional $\psi\colon A\to \mathbb{F}$ and some $A$-differentiable function $f$ defined on an open, simply connected set $U\subset B$. Such $v$ will be called a \emph{component function}. If $B=A$, then such conditions are described in \cite{W1}. We follow essentially the same lines as therein, except the main theorem.

\subsection{Algebraic preliminaries}
Before we come to the main point of this section, we shall now recall some algebraic notions. This part comprises an extension of some facts from \cite{W1}.

\begin{definition}
A commutative algebra $A$ is called a \emph{Frobenius algebra} if there is a linear functional $\phi \colon A \to \mathbb{F}$ such that the bilinear form $(x,y)\mapsto \phi(xy)$ is nondegenerate. 
\end{definition}

Alternatively, one could say that for every linear functional on $A$ there exists unique $y\in A$ such that the functional is of the form $x\mapsto \phi(xy)$. 

We refer the reader to \cite{S} for an extensive account of Frobenius algebras. 

\begin{example}\label{E:Frobenius}
Let $A=\mathbb{F}[x]/(P(x))$ be a quotient algebra of some polynomial $P$. Then $A$ is a Frobenius algebra. 
\end{example}

Next we consider linear and bilinear functions on $A$-modules.

\begin{lemma}\label{L:linear}
Let $A$ be a Frobenius algebra and $B$ be a finitely generated free module over $A$. Let $K\colon B\to \mathbb{F}$ be a linear functional. Then there exists a unique $c\in B$ such that $K(x)=\phi\big(\langle x,c\rangle_A\big)$ for all $x\in B$.
\end{lemma}
\begin{proof}
For any $i=1,\dotsc,n $ the map $x\mapsto K(xb_i)$ is linear, thus it is equal to $\phi(xc_i)$ for some unique $c_i\in A$. Then
\begin{equation*}
K(\sum_{i=1}^nx_ib_i)=\sum_{i=1}^{n}\phi(x_ic_i)=\phi\big(\langle x,c\rangle_A\big).
\end{equation*}
\end{proof}

\begin{lemma}\label{L:bilinear}
Let $A$ and $B$ be as in the preceding lemma. Let $L\colon B\times \dotsc \times B \to \mathbb{F}$ be a multilinear form such that $L(x_1,\dotsc, ax_l,\dotsc,x_r)=L(x_1, \dotsc,ax_k, \dotsc,x_r)$ for all $x_1,\dotsc,x_r \in B$, all $a\in A$ and all $k,l=1,\dotsc,r$. Then there exist unique elements $c_{i_1,\dotsc,i_r}\in A$, $i_1,\dotsc,i_r=1,\dotsc, n$ such that 
\begin{equation*}
L(x^1,\dotsc,x^r)=\phi\left(\sum_{i_1,\dotsc,i_r=1}^n x^1_{i_1}\dotsm x^r_{i_r}c_{i_1,\dotsc,i_r}\right)
\end{equation*}
for all $x^i=\sum_{j=1}^nx^i_jb_j$, $i=1,\dotsc,r$.
\end{lemma}
\begin{proof}
The map $x\mapsto L(xb_{i_1}, b_{i_2},\dotsc,b_{i_r})$ is linear, thus it is equal to $\phi(xc_{i_1,\dotsc,i_r})$ for some unique $c_{i_1,\dotsc,i_r}\in A$. We have
\begin{equation*}
L(a_1b_{i_1},a_2b_{i_2},\dotsc,a_rb_{i_r})=L(a_1\dotsm a_r b_{i_1},  b_{i_2},\dotsc,b_{i_r}).
\end{equation*}
Therefore 
\begin{equation*}
L\bigg(\sum_{i_1=1}^nx^1_{i_1}b_{i_1},\sum_{i=1}^nx^2_{i_2}b_{i_2},\dotsc,\sum_{i=1}^nx^r_{i_r}b_{i_r}\bigg)=\sum_{i_1,\dotsc,i_r=1}^n \phi(x^1_{i_1}\dotsm x^r_{i_r}c_{i_1,\dotsc,i_r}).
\end{equation*}
\end{proof}

\begin{example}
If $r=2$, then $(c_{i_1,i_2})_{i_1,i_2=1,\dotsc,n}$ form a matrix $\tens{C}\in M_{n\times n}(A)$. Then $L$ is given by $L(x,y)=\phi\big(\langle x,\tens{C}y\rangle_A\big)$.
\end{example}

\subsection{Components and generalised Laplace equations}

\begin{lemma}\label{L:components}
Let $v\colon U \to \mathbb{F}$ be a component function of a twice continuously $A$-differentiable map $f\colon U\to A$. Then 
\begin{equation*}
\mathrm{D}^2v(b)(ax,y)=\mathrm{D}^2v(b)(x,ay)
\end{equation*}
for all $x,y\in B$ and $a\in A$.
\end{lemma}
\begin{proof}
By assumption, $v=\psi(f)$ for some linear functional $\psi$. By Lemma \ref{L:secondderivative} we have 
\begin{equation*}
\mathrm{D}^2v(b)(ax,y)=\psi(\mathrm{D}^2f(b)(ax,y))=\psi(\mathrm{D}^2f(b)(x,ay))=\mathrm{D}^2v(b)(x,ay).
\end{equation*}
\end{proof}

\begin{definition}\label{D:short-path}
Let $U\subset B$. We say that $U$ is \emph{short-path connected}, if  there exists a function  $h\colon [0,\infty)\to \mathbb{R}$ such that $h(0)=0$, $h$ is continuous in $0$ and such that for any points $x, y\in U$ there exists a path $\gamma \subset U$ connecting $x$ and $y$ and such that its length $\abs{\gamma}$ satisfies $\abs{\gamma}\leq h(\norm{x-y})$.
\end{definition}

\begin{theorem}\label{T:components}
Suppose that $A$ is a Frobenius algebra and let $B$ be an $A$-module. Assume that $U\subset B$ is an open and simply connected set. Let $t\geq 0$ be a natural number. Let $v\colon U\to \mathbb{F}$ be a $(t+2)$-times continuously differentiable function such that 
\begin{equation*}
\mathrm{D}^2v(b)(ax,y)=\mathrm{D}^2v(b)(x,ay)
\end{equation*}
for all $b\in U$ and $x,y \in B$ and $a\in A$. Then $v=\phi(f)$ for some $(t+2)$-times continuously $A$-differentiable function $f$. Such $f$ is uniquely determined by $v$, up to a constant.
If $U$ is short-path connected and  $v$ is $(t+2)$-times continuously differentiable on $\overline{U}$, then $f$ is also $(t+2)$-times continuously differentiable on $\overline{U}$.
\end{theorem}

\begin{lemma}\label{L:free modules}
The statement of Theorem \ref{T:components} holds true for free modules.
\end{lemma}
\begin{proof}
By Lemma \ref{L:linear} we see that for all $b\in U$ we have $\mathrm{D}v(b)(x)=\phi(\langle x,g(b)\rangle_A)$ for a uniquely determined element $g(b)$. By Lemma \ref{L:bilinear} we have 
\begin{equation*}
\mathrm{D}^2v(b)(x,y)=\phi(\langle x, \tens{C}(b)y\rangle_A)\text{ for a uniquely determined }\tens{C}(b)\in M_{n\times n}(A).
\end{equation*} 
We shall show that $g\colon U\to B$ is $A$-differentiable and that  $\mathrm{D}g(b)(y)=\tens{C}(b)y$. For all $x\in B$ we have
\begin{equation*}
\lim_{y \to 0}  \tfrac{1}{\norm{y}} \left(\mathrm{D}v(b+y)(x)-\mathrm{D}v(b)(x)-\mathrm{D}^2v(b)(x,y)\right)=0.
\end{equation*}
By our previous observations
\begin{equation*}
\lim_{y \to 0}\phi\left(\big\langle x, \tfrac{1}{\norm{y}}(g(b+y)-g(b)-\tens{C}(b)y)\big\rangle_A\right)=0.
\end{equation*}
That means, by Lemma \ref{L:linear}, that the quotients composed with \emph{any} linear functional $\frac1{\norm{y}}(g(b+y)-g(b)-\tens{C}(b)y)$ converge to zero . Since $B$ is finite-dimensional, as a vector space over $\mathbb{F}$, they converge in norm. This is to say, $g$ is differentiable and $\mathrm{D}g(b)(y)=\tens{C}(b)y$. Thus the derivative is $A$-linear.

Observe that $\tens{C}(\cdot)$ is continuous on $U$. Indeed, as $\mathrm{D}^2v(\cdot)(x,y)$ is continuous for all $x,y$, if $b_n\to b$ then 
\begin{equation}\label{eqn:cont}
\phi(\langle x,\tens{C}(b_n)y\rangle_A)=\mathrm{D}^2v(b_n)(x,y)\to \mathrm{D}^2v(b)(x,y)=\phi(\langle x,\tens{C}(b)y\rangle_A).
\end{equation}
Lemma \ref{L:linear} implies that $\tens{C}(b_n)y$ converges to $\tens{C}(b)y$ for all $y$. In consequence $\tens{C}(b_n)$ converges to $\tens{C}(b)$. Hence $g$ is continuously differentiable.

Theorem \ref{T:secondder} tells us that there exists an $A$-differentiable $f\colon U\to A$ such that $\mathrm{D}f(b)(x)=\langle x, g(b)\rangle_A$ for all $b\in U$ and $x\in B$. Such $f$ is unique up to a constant. We have
\begin{equation*}
\mathrm{D}\phi(f(b))(x)=\phi\left(\mathrm{D}f(b)(x)\right)=\phi(\langle x,g(b)\rangle_A)=\mathrm{D}v(b)(x).
\end{equation*}
Thus $\phi(f)$ and $v$ differ only by a constant. Adding a suitable constant to $f$ we get $v=\phi(f)$. 

We have to show that $f$ is $(t+2)$-times continuously differentiable provided that  $v$ is. For all $2\leq k\leq t+2$ the linear form 
\begin{equation*}
(x_1,\dotsc,x_k)\mapsto \mathrm{D}^kv(\cdot)(x_1,\dotsc,x_k)
\end{equation*}
satisfies the assumptions of Lemma \ref{L:bilinear}. Indeed, as $\mathrm{D}^2v(\cdot)$ does, we have\footnote{ We denote $(x_1,\dotsc, x_{l-1},x_{l+1},\dotsc, x_{p-1},x_{p+1},\dotsc, x_k)$ by $(x_1,\dotsc, \hat{x_l},\dotsc, \hat{x_p},\dotsc, x_k)$.}
\begin{equation*}
\begin{aligned}
&\mathrm{D}^kv(\cdot)(x_1,\dotsc, ax_l,\dotsc,x_k)=\mathrm{D}^{k-2}(\mathrm{D}^2v(\cdot)(ax_l,x_p))(x_1,\dotsc, \hat{x_l},\dotsc, \hat{x_p},\dotsc, x_k)=\\
&=\mathrm{D}^{k-2}(\mathrm{D}^2v(\cdot)(x_l,ax_p))(x_1,\dotsc, \hat{x_l},\dotsc, \hat{x_p},\dotsc, x_k)=\mathrm{D}^kv(\cdot)(x_1,\dotsc, ax_p,\dotsc,x_k)
\end{aligned}
\end{equation*}
for all $l<p \leq k$. 
Thus there exist $c_{i_1,\dotsc,i_k}(b)\in A$, $i_1,\dotsc,i_k=1,\dotsc, n$, such that
\begin{equation*}
\mathrm{D}^kv(b)(x^1,\dotsc,x^k)=\phi\left(\sum_{i_1,\dotsc,i_k=1}^n x^1_{i_1}\dotsm x^k_{i_k}c_{i_1,\dotsc,i_k}(b)\right)
\end{equation*}
Assume that we have shown that $f$ is $k$-times continuously differentiable  for some $k<2+t$ and 
\begin{equation}\label{eqn:k-diff}
\mathrm{D}^kf(b)(x^1,\dotsc,x^k)=\sum_{i_1,\dotsc,i_k=1}^n x^1_{i_1}\dotsm x^k_{i_k}c_{i_1,\dotsc,i_k}(b).
\end{equation}
Then, as
\begin{equation*}
\tfrac1{\norm{z}}( \mathrm{D}^kv(b+z)(x^1,\dotsc,x^k)-\mathrm{D}^kv(b)(x^1,\dotsc,x^k)-\mathrm{D}^{k+1}v(b)(x^1,\dotsc,x^k,z))\xrightarrow[{z\to 0}]{}0
\end{equation*}
we have that 
\begin{equation*}
\phi\Bigg(\tfrac1{\norm{z}}\bigg(\sum_{i_1,\dotsc,i_k=1}^n x^1_{i_1}\dotsm x^k_{i_k}\Big(c_{i_1,\dotsc,i_k}(b+z)-c_{i_1,\dotsc,i_k}(b)-\sum_{i_{k+1}=1}^n z_{i_{k+1}}c_{i_1,\dotsc,i_k,i_{k+1}}(b)\Big)\bigg)\Bigg)
\end{equation*}
converges to $0$ as $z$ converges to $0$.
By Lemma \ref{L:linear} we infer that $c_{i_1,\dotsc,i_k}(\cdot)$ are differentiable and
\begin{equation*}
\mathrm{D}c_{i_1,\dotsc,i_k}(\cdot)(z)=\sum_{i=1}^nz_ic_{i_1,\dotsc,i_k,i}(b).
\end{equation*}
Thus the condition (\ref{eqn:k-diff}) is true also for $k+1$. Induction shows that $f$ is $(t+2)$-times differentiable.
Arguing as previously (see (\ref{eqn:cont})) we show that $f$ is $(t+2)$-times continuously differentiable. The same argument shows that if $v$ is $(t+2)$-times continuously differentiable on $\overline{U}$, then all derivatives of $f$, up to order $2+t$, extend continuously to $\overline{U}$. 

The function $f$ extends to $\overline{U}$, since its derivative does. Indeed, let $b_0\in  \overline{U}\setminus U$ and let a sequence $(b_n)_{n=1}^{\infty}\subset U$ converge to $b_0$. Then $(f(b_n))_{n=1}^{\infty}$ is a Cauchy sequence. Indeed, for a curve $\gamma\subset U$ connecting $b_n$ and $b_m$ which satisfies $\abs{\gamma}\leq h(\norm{b_n-b_m})$, we have 
\begin{equation*}
\norm{f(b_m)-f(b_n)}=\Big\lVert\int_0^1 \mathrm{D}f(\gamma(t))(\gamma'(t))\rd t\Big\rVert\leq M \abs{\gamma}\leq M h(\norm{b_n-b_m}),
\end{equation*}
as $\mathrm{D}f(\cdot)$ remains bounded by some $M$ as $n,m\rightarrow \infty$.
Thus we may define 
\begin{equation*}
f(b_0)=\lim_{n\to \infty}f(b_n).
\end{equation*}
This definition does not depend on the choice of the sequence $(b_n)_{n=1}^{\infty}$. Indeed, choose a sequence $(c_n)_{n=1}^{\infty}\subset U$ converging to $b_0$. Then, as before,
\begin{equation*}
\norm{f(c_n)-f(b_n)}\leq M h(\norm{c_n-b_n}).
\end{equation*} 
Thus 
\begin{equation*}
\lim_{n\to \infty}f(b_n)=\lim_{n\to \infty}f(c_n).
\end{equation*}
In particular, for any sequence $(d_n)_{n=1}^{\infty}\subset U$ converging to $b_0$, we have
\begin{equation*}
f(b_0)=\lim_{n\to \infty}f(d_n).
\end{equation*} 
Thus $f$ is continuous on $\overline{U}$.
\end{proof}

\begin{remark}
The assumption about short-path connectedness is used only to show that $f$ itself (and not its derivatives) extends continuously on $\overline{U}$.
\end{remark}

\begin{remark}
We have actually shown that $v=\phi(f)$, where $\phi$ is a linear functional making the algebra Frobenius.
\end{remark}

\begin{proof}[Theorem \ref{T:components}]
Let $v\colon U \to \mathbb{F}$ be as in the theorem, $(t+2)$-times continuously differentiable. Let $\eta\colon A^n\to B$ be a surjective $A$-linear map. Consider the composition
\begin{equation*}
v\circ\eta \colon \eta^{-1}(U)\to \mathbb{F}.
\end{equation*}
Then $v\circ\eta$ satisfies the assumptions of Lemma \ref{L:free modules}. Indeed,
\begin{equation*}
\begin{aligned}
\mathrm{D}^2(v\circ \eta)(b)(ax,y)&=\mathrm{D}(\mathrm{D}(v\circ\eta)(b)(ax))(y)=\mathrm{D}(\mathrm{D}v(\eta(b))(\eta(ax)))(y)=\\
&=\mathrm{D}^2v(\eta (b))(a\eta(x),\eta(y))=\mathrm{D}^2v(\eta (b))(\eta(x),a\eta(y))=\\
&=\mathrm{D}^2(v\circ \eta)(b)(x,ay).
\end{aligned}
\end{equation*}
Thus there exists a $(t+2)$-times continuously $A$-differentiable function $f\colon \eta^{-1}(U)\to A$ such that $v\circ \eta=\phi (f)$.

We shall show that there exists $g\colon U \to A$ such that $g\circ \eta=f$.
Observe that 
\begin{equation*}
\phi(a\mathrm{D}f(b)(x))=\mathrm{D}v(\eta(b))(\eta(ax))=\mathrm{D}v(\eta(b))(a\eta(x))
\end{equation*}
for all $b\in \eta^{-1}(U)$, $x\in A^n$, $a\in A$. Thus if $x\in\ker\eta$, then $\phi(aDf(b)(x))=0$ for all $a\in A$. Since $A$ is Frobenius, $\mathrm{D}f(b)(x)=0$. 

Choose two points $b_1,b_2\in \eta^{-1}(x)$. Then $b_2-b_1\in \ker \eta$. For any $t\in [0,1]$, $b_1+t(b_2-b_1)\in \eta^{-1}(U)$. Thus 
\begin{equation*}
f(b_2)-f(b_1)=\int_0^1 Df(b_1+t(b_2-b_1))(b_2-b_1)\rd t=0
\end{equation*}  
We may thus define $g$ on $U$ by $g(\eta(x))=f(x)$. In view of Lemma \ref{L:etadifferentiability}, we see that $g$ is $(t+2)$-times continuously $A$-differentiable. Moreover 
\begin{equation*}
v\circ \eta=\phi(f)=\phi(g)\circ \eta.
\end{equation*}
As $\eta$ is surjective, we have $v=\phi(g)$.

Assume now that $U$ is short-path connected and that $v$ is $(t+2)$-times continuously differentiable on $\overline{U}$. Then the composition $v\circ\eta$ is again $(t+2)$-times continuously differentiable on $\overline{\eta^{-1}(U)}$
and $\eta^{-1}(U)$ is short-path connected. Therefore $f$ is $(t+2)$-times continuously differentiable on $\overline{\eta^{-1}(U)}$, and by Lemma \ref{L:etadifferentiability}, $g$ is $(t+2)$-times continuously differentiable on $\overline{U}$.
\end{proof}

\begin{example}
Equations
\begin{equation*}
\mathrm{D}^2v(b)(ax,y)=\mathrm{D}^2v(b)(x,ay)
\end{equation*}
are called \emph{generalised Laplace equations}. If $A=\mathbb{C}$, treated as $an \mathbb{R}$-algebra, and $B=A$, they
are equivalent to the Laplace equation in two variables
\begin{equation*}
\frac{\partial^2v}{\partial y^2}(b)=D^2v(b)(i,i)=D^2v(b)(i^2,1)=-\frac{\partial^2v}{\partial x^2}(b).
\end{equation*}
In this case, the theorem tells us that every harmonic function on $U\subset\mathbb{R}^2$ is a real part of some holomorphic function. From this example one readily sees that the assumption about simple connectedness of $U$ can not be dropped. Indeed, $\log\abs{\cdot}$ is harmonic on $\mathbb{C}\setminus\{0\}$, but it is not the real part of a holomorphic function.

If $A=\mathbb{C}$ and $B=\mathbb{C}^n$, then these equations yield the conditions of pluriharmonicity
\begin{equation*}
\begin{aligned}
&\frac{\partial^2v}{\partial x_k\partial x_l}(b)=-D^2v(b)(i^2e_k,e_l)=-D^2v(b)(ie_k,ie_l)=-\frac{\partial^2v}{\partial y_k\partial y_l}(b),\\
&\frac{\partial^2v}{\partial x_k\partial y_l}(b)=-D^2v(b)(e_k,ie_l)=-D^2v(b)(ie_k,e_l)=-\frac{\partial^2v}{\partial y_k\partial x_l}(b).
\end{aligned}
\end{equation*}
\end{example}

\section{$A$-analiticity and further properties}\label{S:analiticity}

\begin{definition}
Let $B$ be a finitely generated $A$-module. We say that a map $L\colon B^i\to A$ is \emph{symmetric} if
\begin{equation*}
L(x_1,\dotsc,x_i)=L(x_{\sigma(1)},\dotsc,x_{\sigma(i)})
\end{equation*}
for any permutation $\sigma\colon \{1,\dotsc,i\}\to\{1,\dotsc,i\}$.
We say that a map $f\colon B\to A$ is $A$-\emph{polynomial} if $f(x)=\sum_{i=0}^{\infty}L_i(x,\dotsc,x)$ for some $A$-multilinear, symmetric functions $L_i\colon B^i\to A$ of which only a finite number is non-zero. The set of such $A$-polynomials is denoted by $B^{*\infty}$.
\end{definition}

For $L_i(x,\dotsc,x)$ we will simply write $L_i(x^i)$.

\begin{definition}
Let $U\subset B$ be an open subset. Let $f\colon U\to A$. We call $f$ an $A$-\emph{analytic} function if for every point in $b_0\in U$ there exists an open neighbourhood $V\subset U$ of $b_0$, such that for $b\in V$
\begin{equation*}
f(b)=\sum_{i=0}^{\infty}L_i((b-b_0)^i),
\end{equation*}
for some symmetric $A$-multilinear $L_i\colon B^i\to A$ such that for any $b\in V$
\begin{equation*}
\sum_{i=0}^{\infty}\norm{L_i((b-b_0)^i)}<\infty.
\end{equation*}
\end{definition}

\begin{example}\label{exa:analytic}
If $B$ is a free module then any $A$-multilinear symmetric mapping $L_i\colon B^i\to A$ is given by
\begin{equation*}
L_i(x^1,\dotsc,x^i)=L_i\Big(\sum_{j=1}^n x^1_{j_1}b_{j_1},\dotsc,\sum_{j_i=1}^nx^i_{j_i}b_{j_i}\Big)=\sum_{j_1,\dotsc,j_i=1}^n x^1_{j_1}\dotsm x^i_{j_i} c_{j_1,\dotsc,j_i},
\end{equation*}
for some $c_{j_1,\dotsc,j_i}\in A$, $j_1,\dotsc,j_l=1,\dotsc,n$ such that $c_{j_1,\dotsc,j_i}=c_{j_{\sigma(1)},\dotsc,j_{\sigma(i)}}$ for any permutation $\sigma$ on $\{1,\dotsc,i\}$.
If $T\colon B\to A$ is a polynomial in coordinate functions $x_k$, with coefficients in $A$, then
\begin{equation*}
T(x)=\sum_{i=1}^l \sum_{0\leq j_1\leq \dots\leq j_i\leq n} x_{j_1}\dotsm x_{j_i}a_{j_1,\dotsc,j_i},
\end{equation*}
for some $a_{j_1,\dotsc,j_i}\in A$, $0\leq j_1\leq\dotsc\leq n$. 
Put
\begin{equation*}
L_i(x^1,\dotsc,x^i)=\sum_{\sigma \in S_i}\frac1{i!}\sum_{0\leq j_1\leq \dots\leq j_i\leq n} x^{\sigma(1)}_{j_1}\dotsm x^{\sigma(i)}_{j_i}a_{j_1,\dotsc,j_i}.
\end{equation*}
Then $L_i$ are $A$-multilinear and symmetric. Moreover
\begin{equation*}
T(x)=\sum_{i=1}^l L_i(x,\dotsc,x).
\end{equation*}
Thus the set of polynomials in coordinate functions $x_k$, with coefficients in $A$, is equal to $B^{*\infty}$. 
\end{example}

\begin{remark}
Consider the set $B^*$ of $A$-linear functions on $B$ with values in $A$. Suppose that $g_1,\dotsc,g_l$ generate $B^*$ over $A$.
Any map $f\colon B \to A$ of the form 
\begin{equation*}
f(x)=\sum_{k_1,\dotsc,k_l=1}^n a_{k_1,\dotsc,k_l}g_1(x)^{k_1}\dotsm g_l(x)^{k_l}
\end{equation*}
is $A$-polynomial.

If $B$ is a free module, then the previous example shows, that any $A$-polynomial map is of this form.
\end{remark}

In \cite{GM} it is proved that, in the case of a free module, $A$-analytic functions on an open set are exactly $A$-differentiable real-analytic functions. 
We prove that this result actually extends to all finitely generated $A$-modules. 

\begin{theorem}[Taylor's theorem]\label{T:Taylor}
Let $f\colon U \to A$ be $(t+1)$-times continuously differentiable function. Let $x_0, x\in U$ be such that $U$ contains the segment $[x_0,x]$. Then 
\begin{equation*}
f(x)=\sum_{k=0}^{t}\frac1{k!}D^kf(x)((x-x_0)^k) +\int_0^1 \frac{(1-s)^t}{(t+1)!}D^{t+1}f(x_0+s(x-x_0))((x-x_0)^{t+1})\rd s.
\end{equation*} 
\end{theorem}

\begin{lemma}[see \cite{KP}]\label{L:analytic}
Let $U\subset \mathbb{F}^t$. Let $f\colon U\to \mathbb{F}$ be a smooth function. The following conditions are equivalent
\begin{enumerate}[\upshape (i)]
\item $f$ is $\mathbb{F}$-analytic,
\item\label{i:boundedder}for any $x\in U$ there exists $r>0$ and $M,C>0$ such that if $\abs{y_i-x_i}<r$ then
\begin{equation*}
\Big|\frac1{k!}\mathrm{D}^kf(y)(z^k)\Big|\leq MC^k r^{-k}\norm{z}^k
\end{equation*}
for all $z\in \mathbb{F}^t$, $k\in\mathbb{N}$.
\end{enumerate}
\end{lemma}

\begin{theorem}\label{T:analytic}
Let $U\subset B$ be an open subset of an $A$-module $B$. Let $f\colon U\to A$. Then $f$ is $A$-differentiable and $\mathbb{F}$-analytic if and only if it is $A$-analytic.
\end{theorem}
\begin{proof}
We infer as in \cite{GM}. Assume that $f$ is $A$-differentiable and $\mathbb{F}$-analytic. Choose $b_0\in U$. By Lemma \ref{L:analytic} applied to each coordinate of $f$, there exists $r>0$ such that if $\norm{b-b_0}<r$ then 
\begin{equation*}
\norm{\mathrm{D}^kf(b)(y^k)}\leq MC^k\norm{y}^k
\end{equation*}
for all $k\in\mathbb{N}$ and some constants $M,C>0$. 
We bound the remainder term in Taylor's formula 
\begin{equation*}
\begin{aligned}
&\norm{f(y)-\sum_{k=0}^{l}\frac1{k!}\mathrm{D}^kf(x)((y-x)^k) }\leq\\
& \leq\int_0^1 \frac{(1-s)^l}{(l+1)!}\norm{\mathrm{D}^{l+1}f(x+s(y-x))((y-x)^{l+1})}\rd s\leq
 \frac1{l+1} M \left(C\norm{y-x}\right)^{l+1}.
\end{aligned}
\end{equation*}
Thus if $\norm{y-x}<\frac1C$, then the series converges. Since $\frac1{k!}\mathrm{D}^kf(x)$ is $A$-multilinear, $f$ is $A$-analytic.

Assume conversely, that $f$ is $A$-analytic. Then $A$-differentiability follows readily from the definition of $A$-analyticity. Moreover, for any linear functional $\psi$ the composition $\psi(f)$ is $\mathbb{F}$-analytic, by Example \ref{exa:analytic}. This means that $f$ itself is $\mathbb{F}$-analytic.
\end{proof}

Let us recall some basics of algebra. We refer the reader to \cite{JAC} for the background.

\begin{definition}
An algebra $A$ is \emph{local} if it has exactly one maximal ideal.
\end{definition}

Suppose that $\mathfrak{m}\subset A$ is a maximal ideal. It is a well-known fact and an easy observation that the quotient $A/\mathfrak{m}$ is a field.

\begin{theorem}\label{T:local}
Let $A$ be a commutative finite dimensional algebra over field $\mathbb{F}$. Then there are finitely many ideals $A_i\subset A$ such that $A=\bigoplus_{i=1}^m A_i$. Each $A_i$ is a local algebra. 
\end{theorem}

Let $A= \bigoplus_{i=1}^{m} A_i$. Then the unit $e\in A$ decomposes into a sum
\begin{equation*}
e=\sum_{i=1}^m e_i,
\end{equation*}
where $e_i\in A_i$. Moreover $e_ie_j=0$ if $i\neq j$. Since $e_ia_i=a_i$ for $a_i\in A_i$, $e_i$ is a unit in $A_i$. Since $ea_i=e_ia_i$, $A_i=e_iA$. 

For any $A$-module $B$ we have $B=\bigoplus_{i=1}^{m}e_iB$. Each $B_i=e_iB$ is an $A_i$-module. Define projection $\pi_{A_i}\colon A\to A$ by $\pi_{A_i}(a)=e_ia$ and $\pi_{B_i}\colon B\to B$ by $\pi_i(b)=e_ib$.

Let $L\colon B\to A$ be an $A$-linear function. Then 
\begin{equation*}
L(b)=\sum_{i=1}^{m}e_iL(b)=\sum_{i=1}^{m}e_iL(e_ib)=\sum_{i=1}^{m}\pi_{A_i}(L(\pi_{B_i}(b))).
\end{equation*}
Functions $\pi_{A_i}\circ L\colon B_i\to A_i$ are $A_i$-linear. 

\begin{proposition}\label{P:spliting}
Assume that $A= \bigoplus_{i=1}^{m} A_i$ is a direct product of algebras. Let $U$ be a convex and open set in a finitely generated $A$-module $B$. Let $f\colon U \to B$ be $A$-differentiable.  Then 
\begin{equation}\label{eqn:form}
f=\sum_{i=1}^{m} f_i\circ\pi_{B_i}
\end{equation}
for some $A_i$-differentiable functions $f_i\colon \pi_{B_i}(U)\to A_i$.
Conversely, any function of this form is $A$-differentiable.
\end{proposition}
\begin{proof}
Let $e\in A$ be the unit. Then, as above, $e=\sum_{i=1}^me_i$ for $e_i\in A_i$ being units of $A_i$. 

Define $g_i=e_if$. If, for some $x\in B$, $e_ix=0$ then
\begin{equation*}
\mathrm{D}g_i(y)(x)=\mathrm{D}(e_if)(y)(x)=\mathrm{D}(e_if)(y)(e_ix)=0,
\end{equation*}
since $e_i^2=e_i$.
By convexity, $g_i(x)=g_i(y)$ if $\pi_{B_i}(x)=\pi_{B_i}(y)$. We define function
\begin{equation*}
f_i\colon \pi_i(U)\to A_i
\end{equation*}
by the formula $f_i(x)=g_i(t)$, where $t\in \pi_i^{-1}(x)$. That is, $f_i\circ \pi_{B_i}=g_i$. Then Lemma \ref{L:etadifferentiability} implies that $f_i$ are $A_i$-differentiable. 
It follows that
\begin{equation*}
f(x)=\sum_{i=1}^{m}e_if(x)=\sum_{i=1}^{m}f_i(\pi_{B_i}(x)).
\end{equation*}

To prove the converse, observe that the derivative of a function of the form (\ref{eqn:form}), is given by 
\begin{equation*}
\mathrm{D}f(u)(x)=\sum_{i=1}^{m} \mathrm{D}f_i(\pi_{B_i}(u)(\pi_{B_i}(x))
\end{equation*}
so 
\begin{equation*}
\mathrm{D}f(u)(ax)=\sum_{i=1}^{m} \mathrm{D}f_i(\pi_{B_i}(u))(e_iax)=\sum_{i=1}^{m} e_ia\mathrm{D}f_i(\pi_{B_i}(u))(e_ix)=a\mathrm{D}f(u)(x).
\end{equation*}
Thus $f$ is $A$-differentiable.
\end{proof}

\begin{lemma}[see \cite{GM}]\label{L:sum}
Let $A$ be a finite dimensional commutative local algebra over $\mathbb{F}\in\{\mathbb{R},\mathbb{C}\}$. Then $A=A/\mathfrak{m}\oplus \mathfrak{m}$ and $\pi_A\colon A\to A/\mathfrak{m}$ is equal to the identity on the first summand and zero on the second.
\end{lemma}

As is shown in \cite{GM}, if there exists a maximal ideal $\mathfrak{m}$ in $A$, such that $A/{\mathfrak{m}} \cong \mathbb{R}$, then there exists a smooth $A$-differentiable function on a free module which is not $A$-analytic. Actually, this result holds for arbitrary modules over Frobenius algebras.

On the other hand, if $A/ \mathfrak{m}\cong \mathbb{C}$ for all maximal ideals, then $A$ is a $\mathbb{C}$-algebra (see Lemma \ref{L:sum}) and from Proposition \ref{P:A-continuous}, one infers that any $A$-continuous function is $A$-analytic. 

\begin{theorem}\label{T:C-continuous}
Let $A$ be a $\mathbb{C}$-algebra. Then every $A$-continuous function is $A$-analytic.
\end{theorem}
\begin{proof}
Let $f\colon U\to B$ be an $A$-continuous function. Then, by Proposition \ref{P:A-continuous}, for any $b\in B$ and any $j=1,\dotsc,n$, there exists $A$-differentiable function $g$ such that 
\begin{equation*}
\mathrm{D}g(y)(x)=xf\circ\eta\circ t^b_j(y).
\end{equation*}
Since $\mathbb{C}\subset A$, the derivative $\mathrm{D}g(y)$ is $\mathbb{C}$-linear, which means that $g$ is $\mathbb{C}$-differentiable and hence $\mathbb{C}$-analytic. Thus $f\circ \eta\circ t^b_j$ is smooth and $A$-differentiable. Then $f\circ \eta$ is smooth and by Proposition \ref{P:separately}, $f\circ \eta$ is $A$-differentiable. Then by Lemma \ref{L:etadifferentiability}, $f$ is $A$-differentiable and since $\mathbb{C}\subset A$ it is $\mathbb{C}$-analytic. By Theorem \ref{T:analytic}, it is $A$-analytic.
\end{proof}

\section{Banach algebras of $A$-differentiable functions}

For any two $A$-differentiable, $A$-valued functions $f,g$, their product $fg$ satisfies
\begin{equation*}
\begin{aligned}
\mathrm{D}(fg)(\cdot)(ax)&=f\mathrm{D}g(\cdot)(ax)+g\mathrm{D}f(\cdot)(ax)=a(f\mathrm{D}g(\cdot)(x)+g\mathrm{D}f(\cdot)(x))=\\
&=a\mathrm{D}(fg)(\cdot)(x).
\end{aligned}
\end{equation*}
Thus, it is again $A$-differentiable. Note the essential role played in the equation by commutativity of $A$.
We see that $A$-differentiable functions form an algebra and it is natural to equip this algebra with a norm $\norm{\cdot}$. We would like the formed algebra to be a Banach algebra. We refer the reader to \cite{Z} for an account on Banach algebras.

\begin{definition}
Assume that $U\subset B$ is open and bounded. Let $C(\overline{U}, A)$ denote the set of all $A$-valued continuous functions on $\overline{U}$ equipped with the Banach algebra norm 
\begin{equation*}
\norm{f}=\sup_{x\in \overline{U}} \norm{f(x)}. 
\end{equation*}
We define the following subsets of $C(\overline{U}, A)$:
\begin{enumerate}[\upshape (i)]
\item 
 $C_A(\overline{U}, A)$ -- the set of all $A$-continuous functions on $U$,
\item $C^{\infty}_A(\overline{U}, A)$ -- the set of all smooth $A$-differentiable functions on $U$, with all derivatives continuous up to the boundary. 
\item $C^{\omega}_A(\overline{U},A)$ -- the set of all $\mathbb{F}$-analytic, $A$-differentiable functions on $U$, with all derivatives continuous up to the boundary. 
\item $C^k_A(\overline{U}, A)$ -- the set of all $k$-times continuously $A$-differentiable functions on $U$, with all derivatives, of order up to $k$, continuous up to the boundary. 
\end{enumerate}
\end{definition}

In the space $C^k_A(\overline{U}, A)$ we introduce a norm given by the formula
\begin{equation*}
\norm{f}_k=\sup_{x\in\overline{U}} \norm{f(x)}+\sum_{i=1}^k\frac1{i!}\sup_{x\in\overline{U}}\sup_{\norm{y}=1} \norm{D^if(x)((y)^i)}.
\end{equation*}
This norm makes $C^k_A(\overline{U}, A)$ a~Banach algebra.

We shall also write $C^0_A(\overline{U},A)$ for $C_A(\overline{U},A)$.

\begin{proposition}\label{P:density}
Assume that $U\subset B$ is an open, convex and bounded set. Then for any natural $k\geq 0$ the algebra $C^{\infty}_A(\overline{U}, A)$ is dense in $C^k_A(\overline{U}, A)$. 
\end{proposition}
\begin{proof}
Let us choose a function $f\colon \overline{U}\to A$ belonging to $C^k_A(\overline{U},A)$. Recall that this means that $f$ is sufficiently differentiable and
\begin{equation*}
\int_{\gamma}(f\circ\eta\circ t^b_j) \rd x=0
\end{equation*}
for all $j=1,\dotsc,n$, all $b\in B$ and all smooth closed curves $\gamma$ in $ (\eta\circ t^b_j)^{-1}(U)$. 

Choose $b_0\in U$ and $0<\delta<1$. Let
\begin{equation*}
U_{\delta}=\{b\in B\colon (1-\delta)(b-b_0)+b_0\in U\}.
\end{equation*}
Then $\delta \leq \mathrm{dist}(B\setminus U_{\delta}, U)$. 
For any $0<\delta<1$ we define $f^{\delta}\colon \overline{U}_{\delta}\to A$ by 
\begin{equation*}
f^{\delta}(b)=f((1-\delta)(b-b_0)+b_0).
\end{equation*}
Then  $f^{\delta}\in C^k_A(\overline{U_{\delta}},A)$ and the sequence $(f^{\delta})_{\delta>0}$ converges in $\norm{\cdot}_k$ to $f$ on $\overline{U}$. Thus it is enough to approximate every $f^{\delta}$ on $\overline{U}$.  

Multiplying $f^{\delta}$ by a smooth bump function equal to one on $U_{\delta/4}$ and zero on $U_{3\delta/4}$ 
we obtain a compactly supported $k$-times continuously diffrentiable function $\tilde{f}^{\delta}$ defined on whole $B$. Again, it is enough to approximate $\tilde{f}^{\delta}$.

Choose a smooth non-negative function $\theta\colon B\to \mathbb{F}$ supported on the unit ball and such that $\int_B\theta\rd\lambda=1$. Here $\lambda$ denotes the Lebesgue measure. Set 
\begin{equation*}
\theta_{\epsilon}(x)={\epsilon}^{t}\theta\Big(\frac{x}{\epsilon}\Big),
\end{equation*}
where $t=\mathrm{dim}_{\mathbb{R}}B$. Then, by uniform continuity of $\tilde{f}^{\delta}$ and its derivatives,  
\begin{equation*}
\tilde{f}^{\delta}_{\epsilon}=\tilde{f}^{\delta}*\theta_{\epsilon}
\end{equation*}
converge in $\norm{\cdot}_k$ on $B$ to $\tilde{f}^{\delta}$ as $\epsilon$ goes to zero. Moreover $\tilde{f}^{\delta}_{\epsilon}$ are smooth. We shall now prove that every $\tilde{f}^{\delta}_{\epsilon}$ is $A$-continuous on $U$ provided that $\epsilon$ is sufficiently small, depending on $\delta$. For this, let $\epsilon<C\delta$, where a constant $C$ depends only on the choice of norms, and choose a smooth closed curve $\gamma\colon [0,1]\to (\eta\circ t^b_j)^{-1}(U)$. 

Then we have by the Fubini's theorem
\begin{equation*}
\begin{aligned}
&\int_{\gamma}(\tilde{f}^{\delta}_{\epsilon}\circ\eta \circ t^b_j)\rd x=\int_0^1 \int_{\mathbb{R}^t}\epsilon^t\theta\Big(\frac{z}{\epsilon}\Big)\tilde{f}^{\delta}(\eta( t^b_j(\gamma(t)))-z)\rd\lambda(z) \rd\gamma(t)=\\
&=\int_{B(0,\epsilon)}\epsilon^t\theta\Big(\frac{z}{\epsilon}\Big)\int_0^1 \tilde{f}^{\delta}(\eta(t^b_j(\gamma(t)))-z)\rd\gamma(t)\rd\lambda(z)=\\
&=\int_{B(0,\epsilon)}\epsilon^t\theta\Big(\frac{z}{\epsilon}\Big)\int_{\gamma} (f^{z,\delta}\circ\eta\circ t^b_j)\rd x\rd\lambda
\end{aligned}
\end{equation*} 
where $f^{z,\delta}$, given by the formula $f^{z,\delta}(a)=f^{\delta}(a-z)$, is an $A$-continuous function,
as $z\in B(0,\epsilon)$. Thus the integral vanishes.
By Corollary \ref{C:equivalence} we see that every $\tilde{f}^{\delta}_{\epsilon}$ is $A$-differentiable.
\end{proof}

\begin{corollary}
If $A$ is a $\mathbb{C}$-algebra, then $C^{\omega}_A(\overline{U},A)$ is a complex Banach algebra.
\end{corollary}
\begin{proof}
By Theorem \ref{T:C-continuous}, any $A$-continuous function is $A$-analytic. Therefore we have $C^{\omega}_A(\overline{U},A)=C_A(\overline{U},A)$. 
\end{proof}

\begin{definition}
Assume that $U\subset A$ is an open set. We say that a function $f\colon U\to A$ admits an $A$-primitive if there exists $g\colon U\to A$ such that $\mathrm{D}g(\cdot)(x)=xf(\cdot)$.
\end{definition}

\begin{corollary}
Assume that $U\subset A$ is an open and simply connected set. Then the set of continuous functions which admit an $A$-primitive forms an algebra.
\end{corollary}
\begin{proof}
Proposition \ref{P:A-continuous} shows that a continuous function admits an $A$-primitive if and only if it is $A$-continuous. The set $C_A(\overline{U},A)$ of $A$-continuous functions is an algebra, since it is the closure of the set of smooth $A$-differentiable functions.
\end{proof}

It is interesting to compare these two situations: the first when $A$ is a $\mathbb{C}$-algebra, and the second when $A$ is a $\mathbb{R}$-algebra. In the first case every $A$-differentiable function is $A$-analytic. It may be differentiated as many times as desired, as well primitives may be taken as many times as we want to. In the second case, taking primitives is always possible.

\section{Components of $A$-differentiable functions revisited}\label{S:structure}

We shall now provide the precise formulae for $A$-differentiable functions on finitely generated modules over commutative, finite-dimensional algebras.

\subsection{Preliminaries}

Before we concentrate on the introduced Banach algebras, let us present some preliminary lemmas.

\begin{proposition}[see \cite{MO}]
For any $x_1,\dotsc,x_k\in A$ 
\begin{equation*}
x_1\dotsm x_k=\frac1{2^kk!}\sum_{\epsilon_1,\dotsc,\epsilon_k= \pm 1}\Bigg( \epsilon_1\dotsm\epsilon_k \Big(\sum_{i=1}^k \epsilon_i x_i\Big)^k\Bigg).
\end{equation*}
\end{proposition}

\begin{corollary}\label{C:powers}
Let $I\subset A$ be an ideal and let $a\in A$. Assume that $I^ka \neq0$. Then there exists $v\in I$ such that $v^ka\neq 0$.
\end{corollary}
\begin{proof}
Follows from the preceding proposition.
\end{proof}

\begin{lemma}\label{L:polynomials}
Let $U\subset \mathbb{F}^n$ and let $T\subset\mathbb{F}^m$ be an infinite, bounded set. Let $k, k_1,\dotsc,k_m\in\mathbb{N}$. Let
\begin{equation*}
\begin{aligned}
X_{k_1,\dotsc,k_m}=&\Big\{f\colon \overline{U}\times T\to A \colon f(u,t_1,\dotsc,t_m)=\sum_{i_1=0}^{k_1}\dots\sum_{i_m=0}^{k_m}f_{i_1\dotsc i_n}(u)t_1^{i_1}\dotsm t_m^{i_m},\\
& f_{i_1\dotsc i_n}\in C^k(\overline{U},A) \text{ for all }i_1,\dotsc,i_m\Big\}.
\end{aligned}
\end{equation*}
Assume that $(f^p)_{p=0}^{\infty}$ is a sequence in $X_{k_1,\dotsc,k_m}$ such that
\begin{equation*}
\sup_{t\in T}\norm{f^p(\cdot,t)-f(\cdot,t)}_k\to 0.
\end{equation*}
Then for any  $j=1,\dotsc,m$ and any $i_j=0, \dotsc, k_j$ 
\begin{equation*}
\norm{f_{i_1\dotsc i_m}^p-f_{i_1\dotsc i_m}}_k\to 0.
\end{equation*}
\end{lemma}
\begin{proof}
The space 
\begin{equation*}
Y_{k_1,\dotsc,k_m}=\{f\colon T\to A \colon f(t_1,\dotsc,t_m)=\sum_{i_1=0}^{k_1}\dots\sum_{i_m=0}^{k_m}f_{i_1\dotsc i_n}t_1^{i_1}\dotsm t_m^{i_m}\}
\end{equation*}
is finite-dimensional. Thus any two norms are equivalent on $Y_{k_1,\dotsc,k_m}$. There exists a constant $G>0$ such that for any $h\in Y_{k_1,\dotsc,k_m}$
\begin{equation*}
\sup_{ t\in T}\norm{h(t)}\geq G \sum_{i_1=0}^{k_1}\dots\sum_{i_m=0}^{k_m}\norm{h_{i_1\dotsc i_m}}.
\end{equation*}
Thus for any $g\in X_{k_1,\dotsc,k_m}$ and any indices $i_1,\dotsc,i_m$
\begin{equation*}
\sup_{t\in T}\norm{g(\cdot,t)}_k=\sup_{t\in T}\sum_{i=0}^k \frac1{i!}\sup_{u\in \overline{U}}\sup_{\norm{y}=1}\norm{\mathrm{D}^ig(u,t)((y)^i)}
\geq G \norm{g_{i_1\dotsc i_m}}_k.
\end{equation*}
Then
\begin{equation*}
\sup_{ t\in T}\norm{f^p(\cdot,t)-f(\cdot,t)}_k\geq G \norm{f^p_{i_1\dotsc i_m}-f_{i_1\dotsc i_m}}_k.
\end{equation*}
\end{proof}

\subsection{Structure}

By Theorem \ref{T:local} and Proposition \ref{P:spliting} we may now solely consider the case of a local algebra. Let $A$ be a local algebra with the maximal ideal $\mathfrak{m}$, so that $A/\mathfrak{m}$ is a field. For any ideal $I\subset A$ let
\begin{equation*}
\pi_I\colon I\to I/\mathfrak{m}I
\end{equation*}
denote the quotient map. Let $\rho_I\colon  I/\mathfrak{m}I\to I$ be a fixed $A/\mathfrak{m}$-linear map such that $\pi_I\circ \rho_I=\mathrm{id}$. 

Assume that an $A$-module $B$ has a decomposition into two submodules such that $B=C\oplus D$. Let 
\begin{equation*}
\pi_D\colon B\to D/\mathfrak{m}D
\end{equation*}
be an $A$-linear map which is a composition of a projection onto $D$ along $C$ and the quotient map. 
Let 
\begin{equation*}
\rho_D\colon D/\mathfrak{m}D\to B
\end{equation*}
be $A/\mathfrak{m}$-linear map which is a composition of a fixed linear injection from $D/\mathfrak{m}D$ to $D$ and the inclusion of $D$ into $B$ such that $\pi_D\circ\rho_D=\mathrm{id}$. 

\begin{definition}\label{D:L}
We say that a module $B$ and an ideal $I\subset A$ have the \emph{lifting property} if the following conditions are satisfied:
\begin{enumerate}[\upshape (i)]
\item \label{i:comp}there exist submodules $C, D\subset B$ such that $B=C\oplus D$ and such that a symmetric $A/\mathfrak{m}$-multilinear function $f\colon B\times\dotsm\times B \to I/\mathfrak{m}I$ admits a lifting to a symmetric $A$-multilinear homomorphism $h\colon B\times\dotsm \times B\to I$ such that $\pi_I\circ h =f$ and for any $b_1,\dotsc,b_{i-1}\in B$, 
\begin{equation*}
h(b_1,\dotsc,b_{j-1},\cdot, b_{j+1},\dotsc,b_{i-1})|_{C}=0,
\end{equation*} if and only if $f$ is  \emph{admissible}, that is for any $b_1,\dotsc,b_i\in B$, 
\begin{equation*}
f(b_1,\dotsc,b_{j-1},\cdot, b_{j+1},\dotsc,b_i)|_{C}=0;
\end{equation*}
such $h$ we shall call a \emph{lifting} of $f$,
\item \label{i:comp1} for any $j\in\mathbb{N}$, there exists an $A/\mathfrak{m}$-linear operator $\mathrm{G}^j_I$ which assigns to an admissible function $f\colon B\times\dotsm\times B \to I/\mathfrak{m}I$ its lifting 
\begin{equation*}
\mathrm{G}^j_I f \colon B\times\dotsm\times B\to I,
\end{equation*}
\item \label{i:comp2}for any $f\colon D/\mathfrak{m}D\to I/\mathfrak{m}I$ and $v\in B$ we have
\begin{equation*}
\mathrm{G}^1_I(f\circ \pi_D)(v)=\rho_I(f(\pi_D(v))),
\end{equation*}
\item \label{i:comp3} for any admissible $f\colon B\times\dotsm\times B\to I/\mathfrak{m}I$ and any $v_1,\dotsc,v_{j+1}\in B$
\begin{equation*}
\big(\mathrm{G}^j_If(\cdot,\dotsc,\cdot,v_{j+1})\big)(v_1,\dotsc,v_j)=\big(\mathrm{G}^{j+1}_If\big)(v_1,\dotsc,v_j,v_{j+1})
\end{equation*}
\end{enumerate}
\end{definition}

\begin{lemma}\label{L:propertyL}
Let $I\subset A$ be an ideal. Let $i_1,\dotsc,i_k\in I$ give in the quotient an $A/\mathfrak{m}$-basis of $I/\mathfrak{m}I$. Then an $A$-module $B$ and ideal $I$ have the lifting property with decomposition $B=C\oplus D$ if and only if for any $j=1,\dotsc,k$ module the $B$ and ideal $(i_j)$ generated by $i_j$ have the lifting property with decomposition $B=C\oplus D$.
\end{lemma}

\begin{proof}
Let $i_1,\dotsc,i_k\in I$ be as in the statement of the lemma, that is, their images by the quotient map give an $A/\mathfrak{m}$-basis of $I/\mathfrak{m}I$. 
Choose the maps $\rho_{(i_l)}$ in such a way that for each $l=1,\dotsc,k$ there is $\rho_I([i_l])=\rho_{(i_l)}([i_l])$.

Suppose that for each $l$, the module $B$ and the ideal $(i_l)$ have the lifting property, with $A/\mathfrak{m}$-linear operators $G^j_{(i_l)}$ assigning to an admissible function $f_l\colon B\times \dotsc B\to (i_l)/\mathfrak{m}{i_l}$ its lifting $G^j_{(i_l)}f_l$. We would like to define a lifting operator for admissible functions $f\colon B\times \dotsc B\to I/\mathfrak{m}I$.

Any symmetric, admissible, $A/\mathfrak{m}$-multilinear map 
\begin{equation*}
f\colon B\times \dotsc\times B\to I/\mathfrak{m}I
\end{equation*}
can be uniquely written as
\begin{equation*}
f=\sum_{l=1}^k f_l[i_l],
\end{equation*}
where each $f_l[i_l]$ is a symmetric, admissible, $A/\mathfrak{m}$-multilinear map.
Define
\begin{equation*}
G^j_If=\sum_{l=1}^k G^j_{(i_l)}\big(f_l[i_l]\big).
\end{equation*}
Then clearly $G^j_I$ is $A/\mathfrak{m}$-linear, and satisfies properties (\ref{i:comp2}) and (\ref{i:comp3}) of the lemma, as each of the operators $G^j_{(i_l)}$ does satisfy them. Indeed, let us check (\ref{i:comp2}). Take any $f\colon D/\mathfrak{m}D\to I/\mathfrak{m}I$ and $v\in B$. Then
\begin{align*}
&G^j_I(f\circ \pi_D)(v)=\sum_{l=1}^k G^j_{(i_l)}\big(f_l\circ \pi_D[i_l]\big)(v)=\sum_{l=1}^k \rho_{(i_l)}(f_l\circ \pi_D[i_l](v))=\\
&=\sum_{l=1}^kf_l\circ \pi_D(v) \rho_{(i_l)}([i_l])=\sum_{l=1}^kf_l\circ \pi_D(v) \rho_{I}([i_l])=\\
&=\rho_I\Big(\sum_{l=1}^k f_l\circ\pi_D(v)[i_l]\Big)=\rho_l (f\circ \pi_D(v)).
\end{align*}

Conversely, suppose that $B$ and $I$ have the lifting property with decomposition $B=C\oplus D$. Let $G^j_I$ be the corresponding lifting operator. Let $l=1,\dotsc, k$ and let $f\colon B\times\dotsc\times B\to (i_l)/\mathfrak{m}(i_l)$ be a symmetric, $A/\mathfrak{m}$-linear, admissible function. Define a map $\eta_l\colon (i_l)/\mathfrak{m}(i_l)\to I/\mathfrak{m}I$ by the formula $\eta_l([ai_l])=[ai_l]$. Clearly, it is well defined and $A/\mathfrak{m}$-linear.
Observe that $f$ is admissible, symmetric and $A/\mathfrak{m}$-multilinear if and only if $\eta_l\circ f$ is. Set $G^j_{(i_l)}(f)=G^j_I(\eta_l\circ f)$. Employing the condition $\rho_I([i_l])=\rho_{(i_l)}([i_l]$ one readily verifies (\ref{i:comp2}). Property (\ref{i:comp3}) is verified readily.
\end{proof}

\begin{lemma}\label{L:passingdown}
Let $B$ be a finitely generated $A$-module such that $B$ and an ideal $I$ have the lifting property. Let $U\subset B$ be an open, convex and bounded set. Assume that $f\colon \overline{U}\to I$ is an $A$-continuous function. Then there exists a unique $A/\mathfrak{m}$-continuous function $g\colon \overline{\pi_D(U)}\to I/\mathfrak{m}I$ such that 
\begin{equation}\label{eqn:interlacing}
\pi_I\circ f=g\circ \pi_D.
\end{equation}
\end{lemma}
\begin{proof}
Assume first that $f$ is $A$-differentiable. 

The map $\pi_I\colon I\to I/\mathfrak{m}I$ is $A$-linear, hence $A$-differentiable. Thus the composition $\pi_I\circ f$ is $A$-differentiable. Let $u\in U$, $b\in B$ and $a\in \mathfrak{m}$. Since we have 
\begin{equation*}
\mathrm{D}(\pi_I\circ f)(u)(b)\in I/\mathfrak{m}I,
\end{equation*}
then
\begin{equation*}
\mathrm{D}(\pi_I\circ f)(u)(ab)=a\mathrm{D}(\pi_I\circ f)(u)(b)=0.
\end{equation*}
Moreover, since $\mathrm{D}(\pi_I\circ f)(u)=\pi_I\circ \mathrm{D}f(u)$, the function $\mathrm{D}(\pi_I\circ f)(u)$ admits a lifting. By the lifting property, 
\begin{equation*}
\mathrm{D}(\pi_I\circ f)(u)|_C=0.
\end{equation*}
Let $u_1,u_2\in U$ be such that $\pi_D(u_1)=\pi_D(u_2)$. By convexity of $U$
\begin{equation*}
\pi_I\circ f(u_1)=\pi_I\circ f(u_2).
\end{equation*}
Therefore, there is $g\colon \pi_D(U)\to I/\mathfrak{m}I$ such that $\pi_I\circ f=g\circ \pi_D$. $A$-differentiability of $g$ follows from Lemma \ref{L:etadifferentiability}, as $\pi_D$ is an $A$-linear surjection. In particular, $g$ is $A/\mathfrak{m}$-differentiable.

If $f$ is merely $A$-continuous, then by Proposition \ref{P:density} we may choose a sequence of smooth $A$-differentiable functions $(f_n)_{n=0}^{\infty}$ converging to $f$ uniformly. Then, by the formula (\ref{eqn:interlacing}), the corresponding functions $(g_n)_{n=0}^{\infty}$ also converge. For its limit (\ref{eqn:interlacing}) still holds. Uniqueness is a consequence of (\ref{eqn:interlacing}).
\end{proof}

Observe that the assignment $f\mapsto g$ is $A/\mathfrak{m}$-linear. Using the lemma we may define the continuous linear operator
\begin{equation*}
\mathrm{H}_I\colon C_A(\overline{U}, I)\to  C_{A/\mathfrak{m}}(\overline{\pi_D(U)}, I/\mathfrak{m}I),
\end{equation*}
by the formula $\mathrm{H}_I(f)=g$.

\begin{definition}
Let $a\in A$. We define a natural number $k(a)$ to be the largest $k\in \mathbb{N}$ such that $\mathfrak{m}^ka\neq \{0\}$.
\end{definition}

\begin{lemma}[see \cite{GM}]\label{L:generate}
There exists a system of elements $e_1,\dotsc,e_r\in\mathfrak{m}$ and a set of multi-indices $M\subset \mathbb{N}^r$, $M\ni (0,\dotsc,0)$, such that:
\begin{enumerate}[\upshape (i)]
\item the cosets of $e^{\alpha}=e_1^{\alpha_1}\dotsm e_r^{\alpha_r}$, $\abs{\alpha}=k$, $\alpha\in M$, modulo $\mathfrak{m}^{k+1}$ are a basis of the $A/\mathfrak{m}$-vector space $\mathfrak{m}^k/\mathfrak{m}^{k+1}$ for every $k\leq l-1$,
\item the elements $(e^{\alpha})_{\alpha \in M}$ are a basis of $A$ as an $A/\mathfrak{m}$-vector space.
\end{enumerate}
\end{lemma}

Let $l\in\mathbb{N}$ be such that $\mathfrak{m}^{l+1}=\{0\}$, but $\mathfrak{m}^l\neq\{0\}$. That is $l=k(1)$.

We choose $\rho_{\mathfrak{m}^k}$ so that $\rho_{\mathfrak{m}^k}([e^{\alpha}])=e^{\alpha}$ for any $\alpha \in M$ such that $\abs{\alpha}=k$.

\begin{remark}
Let $f\in C_A(\overline{U}, \mathfrak{m}^k)$. Observe that\footnote{Here $[a]$ denotes the coset of an element $a\in A$.} $H\mathrm{H}_{\mathfrak{m}^k}f=\sum_{\alpha\in M,\abs{\alpha}=k}f_{\alpha}[e^{\alpha}]$ for some functions $f_{\alpha}$, which we shall, abusing notation, call $\mathrm{H}_{(e^{\alpha})}f$. This shall not lead to any misunderstanding. For if $\abs{\alpha}=k$ and $f\in C_A(\overline{U}, (e^{\alpha}))$, then $\mathrm{H}_{\mathfrak{m}^k}f=\mathrm{H}_{(e^{\alpha})}f$.
\end{remark}

\begin{lemma}\label{L:smoothing}
Let $k,p\in\mathbb{N}$. Let $B$ be a finitely generated $A$-module such that $B$ and $\mathfrak{m}^k$ have the lifting property. Let $U\subset B$ be an open, convex and bounded set. Assume that $f\in C^p_A(\overline{U}, \mathfrak{m}^k)$. Then 
\begin{equation*}
\mathrm{H}_{\mathfrak{m}^k}f\in  \bigoplus_{\mathclap{\alpha\in M, \abs{\alpha}=k}}C^{p+k(e^{\alpha})}_{A/\mathfrak{m}}(\overline{\pi_D(U)}, (e^{\alpha})/\mathfrak{m}(e^{\alpha})).
\end{equation*}
Moreover, the linear operator 
\begin{equation*}
\mathrm{H}_{\mathfrak{m}^k}\colon C^p_A(\overline{U}, \mathfrak{m}^k)\to \bigoplus_{\mathclap{\alpha\in M, \abs{\alpha}=k}} C^{p+k(e^{\alpha})}_{A/\mathfrak{m}}(\overline{\pi_D(U)}, (e^{\alpha})/\mathfrak{m}(e^{\alpha})),
\end{equation*}
is continuous.
\end{lemma}
\begin{proof}
Let $v\in \mathfrak{m}$.
Let $f\in C^p_A(\overline{U}, \mathfrak{m}^k)$ be smooth in $U$. Then for any $u\in U$ and any sufficiently small $b\in B$ by Taylor's theorem we have
\begin{equation*}
f(u+vb)=f(u)+\sum_{j=1}^{k(a)}\frac1{j!}\mathrm{D}^jf(u)((vb)^j)=f(u)+\sum_{j=1}^{l-k}\frac1{j!}v^j\mathrm{D}^jf(u)((b)^j),
\end{equation*}
as by Lemma \ref{L:secondderivative} the derivatives are $A$-multilinear. Observe that $f=\rho_{\mathfrak{m}^k} \pi_{\mathfrak{m}^k}+h$, where $h$ has values in $\mathfrak{m}^{k+1}$. Indeed, $\pi_{\mathfrak{m}^k}(f-\rho_{\mathfrak{m}^k} \pi_{\mathfrak{m}^k}f)=0$.
Then for any $j=1,\dotsc,l-k$
\begin{equation*}
v^j\mathrm{D}^jf(u)((b)^j)=v^j\mathrm{D}^j(\rho_{\mathfrak{m}^k}\pi_{\mathfrak{m}^k} f)(u)((b)^j)+v^j\mathrm{D}^jh(u)((b)^j).
\end{equation*}
By the definition of $\mathrm{H}_{\mathfrak{m}^k}$, $\pi_{\mathfrak{m}^k}\circ f=\mathrm{H}_{\mathfrak{m}^k}f\circ \pi_D$. Thus
\begin{equation*}
v^j\mathrm{D}^jf(u)((b)^j)=v^j\rho_{\mathfrak{m}^k}\mathrm{D}^j(\mathrm{H}_{\mathfrak{m}^k}f)(\pi_D(u))((\pi_D(b))^j)+v^j\mathrm{D}^jh(u)((b)^j).
\end{equation*}
Then
\begin{equation}\label{eqn:poly}
\begin{aligned}
v^j\mathrm{D}^jf(u)((b)^j)&=\sum_{\mathclap{\alpha\in M, \abs{\alpha}=k}}v^je^{\alpha}\mathrm{D}^j(\mathrm{H}_{(e^{\alpha})}f)(\pi_D(u))((\pi_D(b))^j)+\\
&+v^j\mathrm{D}^jh(u)((b)^j).
\end{aligned}
\end{equation}

Choose now a function $g\in C^p_A(\overline{U},\mathfrak{m}^k)$ and a sequence $(g_n)_{n=1}^{\infty}\in C^{\infty}_A(\overline{U},\mathfrak{m}^k)$ which converges to $g$ in the norm $\norm{\cdot}_p$. Then, by Lemma \ref{L:polynomials} for any small $b\in B$ the corresponding left-hand sides of the equation (\ref{eqn:poly}) converge on $\overline{U}$ in the norm $\norm{\cdot}_p$. So do the right-hand sides.  

Therefore for any $v\in\mathfrak{m}$, 
\begin{equation*}
v^je^{\alpha}\mathrm{D}^j(\mathrm{H}_{(e^{\alpha})}f)(\pi_D(u))((\pi_D(b))^j)
\end{equation*}
converge. By Corollary \ref{C:powers} there is $v\in\mathfrak{m}$ such that $v^{k(e^{\alpha})}e^{\alpha}\neq 0$. It follows that $(\mathrm{H}_{(e^{\alpha})}g_n)_{n=1}^{\infty}$ converges in $\norm{\cdot}_{p+k(e^{\alpha})}$. Since $(\mathrm{H}_{(e^{\alpha})}g_n)_{n=0}^{\infty}$ converges to $\mathrm{H}_{(e^{\alpha})}g$ uniformly we have that
\begin{equation*}\mathrm{H}_{(e^{\alpha})}g\in C^{p+k(e^{\alpha})}_{A/\mathfrak{m}}(\overline{\pi_D(U)},(e^{\alpha})/\mathfrak{m}(e^{\alpha})).
\end{equation*}
\end{proof}

\begin{theorem}\label{T:lifting}
Let $p\in\mathbb{N}$. Let $a\in A$. Let $B$ be a finitely generated $A$-module such that $B$ and $(a)$ have the lifting property. Let $U\subset B$ be an open, convex and bounded set. There exists a continuous map 
\begin{equation*}
\mathrm{T}_{(a)}\colon  C^{p+k(a)}_{A/\mathfrak{m}}(\overline{\pi_D(U)}, (a)/\mathfrak{m}(a)) \to C^p_A(\overline{U}, (a))
\end{equation*}
such that $\mathrm{H}_{(a)}\circ \mathrm{T}_{(a)}=\mathrm{id}$.

Moreover, if $k(a)=0$, then $\mathrm{T}_{(a)}\circ \mathrm{H}_{(a)}=\mathrm{id}$. If $b\in A$ is such that $(a)\supset (b)$, $k(a)>k(b)$ and $B$ and $(b)$ have the lifting property, then $\mathrm{H}_{(a)}\circ \mathrm{T}_{(b)}=0$.
\end{theorem}
\begin{proof}
We want to define an operator
\begin{equation*}
\mathrm{T}_{(a)}\colon  C^{p+k(a)}_{A/\mathfrak{m}}(\overline{\pi_D(U)}, (a)/\mathfrak{m}(a)) \to C^p_A(\overline{U}, (a))
\end{equation*}
such that $\mathrm{H}_{(a)}\circ \mathrm{T}_{(a)}=\mathrm{id}$.

Let $g\in C^{p+k(a)}_{A/\mathfrak{m}}(\overline{\pi_D(U)}, (a)/\mathfrak{m}(a))$.
For any $i=1,\dotsc,k(a)$, and any element $p\in \overline{\pi_D(U)}$, the function 
\begin{equation*}
\mathrm{D}^ig(p)\colon \pi_D(B)^i\to (a)/\mathfrak{m}(a)
\end{equation*}
is multilinear and so is its composition with $\pi_D$, 
\begin{equation*}
\mathrm{D}^ig(p)(\pi_D(\cdot),\dotsc,\pi_D(\cdot))\colon B^i\to (a)/\mathfrak{m}(a).
\end{equation*}
For any fixed $b_1,\dotsc,b_{i-1}$, the function $\mathrm{D}^ig(p)(\pi_D(b_1),\dotsc,\pi_D(\cdot),\dotsc,\pi_D(b_{i-1}))$ is linear and vanishes on $C$. By the lifting property there exist liftings 
\begin{equation*}
\mathrm{G}^i_{(a)}\mathrm{D}^ig(p)\colon B^i\to (a).
\end{equation*} 
Define 
\begin{equation}\label{eqn:lift}
\mathrm{T}_{(a)}g(u)=\rho_{(a)}g(\pi_D(u))+\sum_{j=1}^{k(a)}\frac1{j!}\mathrm{G}^j_{(a)}\mathrm{D}^jg(\pi_D(u))((u-\rho_D\pi_D(u))^j).
\end{equation}

Since $g\in C^{p+k(a)}_{A/\mathfrak{m}}(\overline{\pi_D(U)}, (a)/\mathfrak{m}(a))$ we see that $\mathrm{T}_{(a)}g\colon \overline{U}\to (a)$ belongs to $C^p_{A}(\overline{U},A)$. Moreover $\mathrm{T}_{(a)}$ is continuous.

We need to check that $\mathrm{T}_{(a)}g$ is $A$-continuous. For this, let us assume first that $g$ is smooth in $\pi_D(U)$. We shall show that $\mathrm{T}_{(a)}g$ is $A$-differentiable. 

As $g$ is smooth in $\pi_D(U)$ it makes sense to define $\mathrm{G}^i_{(a)}\mathrm{D}^ig$, for all $i\in\mathbb{N}$, in the same manner as before.

Assume first that $v\in B$ is such $\pi_D(v)=0$. Then by multilinearity of $\mathrm{G}^j_{(a)}D^jg$ we see that
\begin{equation*}
\mathrm{D}\mathrm{T}_{(a)}g(u)(v)=\sum_{j=1}^{k(a)}\frac1{(j-1)!}\mathrm{G}^j_{(a)}\mathrm{D}^jg(\pi_D(u))((u-\rho_D\pi_D(u))^{j-1},v).
\end{equation*}
When $v\in B$ is such $\rho_D\pi_D(v)=v$, then
\begin{equation*}
\begin{aligned}
\mathrm{D}\mathrm{T}_{(a)}g(u)(v)=&\rho_{(a)} \mathrm{D}g(\pi_D(u))(\pi_D(v))+\\
&+\sum_{j=1}^{k(a)}\frac1{j!}\mathrm{G}^j_{(a)}\mathrm{D}^{j+1}g(\pi_D(u))(\cdot,\pi_D(v))((u-\rho_D\pi_D(u))^j).
\end{aligned}
\end{equation*}
By the compatibility conditions (\ref{i:comp1})-(\ref{i:comp3}) of the lifting property we see that
\begin{equation*}
\begin{aligned}
\mathrm{D}\mathrm{T}_{(a)}g(u)(v)=&\mathrm{G}^1_{(a)} \mathrm{D}g(\pi_D(u))(v)+\\
&+\sum_{j=1}^{k(a)}\frac1{j!}\mathrm{G}^{j+1}_{(a)}\mathrm{D}^{j+1}g(\pi_D(u))((u-\rho_D\pi_D(u))^j,v).
\end{aligned}
\end{equation*}
Observe that $v=v-\rho_D\pi_D(v)+\rho_D\pi_D(v)$. We have 
\begin{equation*}
\rho_D\pi_D(\rho_D\pi_D(v))=\rho_D\pi_D(v)\text{ and }\rho_D\pi_D(v-\rho_D\pi_D(v))=0.
\end{equation*}
Therefore the derivative of $\mathrm{T}_{(a)}g$ for an arbitrary $v$ is equal to
\begin{equation*}
\begin{aligned}
\mathrm{D}\mathrm{T}_{(a)}g(u)(v)&=\sum_{j=1}^{k(a)}\frac1{(j-1)!}\mathrm{G}^j_{(a)}\mathrm{D}^jg(\pi_D(u))((u-\rho_D\pi_D(u))^{j-1},v-\rho_D\pi_D(v))+\\
&+\sum_{j=0}^{k(a)}\frac1{j!}\mathrm{G}^{j+1}_{(a)}\mathrm{D}^{j+1}g(\pi_D(u))((u-\rho_D\pi_D(u))^j,\rho_D\pi_D(v))=\\
&=\sum_{j=0}^{k(a)}\frac1{j!}\mathrm{G}^{j+1}_{(a)}\mathrm{D}^{j+1}g(\pi_D(u))((u-\rho_D\pi_D(u))^j,v)-\\
&-\frac1{(k(a))!}\mathrm{G}^{k(a)+1}_{(a)}\mathrm{D}^{k(a)+1}g(\pi_D(u))((u-\rho_D\pi_D(u))^j,v-\rho_D\pi_D(v))
\end{aligned}
\end{equation*}
From the definition of $k(a)$ and $A$-multilinearity of $\mathrm{G}^{k(a)+1}_{(a)}D^{k(a)+1}g(\pi_D(u))$ it follows that 
\begin{equation*}
\mathrm{G}^{k(a)+1}_{(a)}\mathrm{D}^{k(a)+1}g(\pi_D(u))((u-\rho_D\pi_D(u))^{k(a)},v-\rho_D\pi_D(v))=0.
\end{equation*}
Thus the derivative of $\mathrm{G}_{(a)}g$ is $A$-linear and, in consequence, $\mathrm{G}_{(a)}g$ is $A$-differentiable.

If $g\in C^{p+k(a)}_{A/\mathfrak{m}}(\overline{\pi_D(U)}, A/\mathfrak{m})$ is not smooth in $\pi_D(U)$, using Proposition \ref{P:density} we choose a sequence 
\begin{equation*}
(g_n)_{n=1}^{\infty}\subset C^{p+k(a)}_{A/\mathfrak{m}}(\overline{\pi_D(U)}, A/\mathfrak{m})
\end{equation*} 
of smooth functions converging in $\norm{\cdot}_{p+k(a)}$ to $g$. Then $\mathrm{T}_{(a)}g_n$ converges to $\mathrm{T}_{(a)}g$ in $C^p(\overline{U},A)$. Since all $\mathrm{T}_{(a)}g_n$ are $A$-continuous, so is $\mathrm{T}_{(a)}g$.

We will now show that $\mathrm{H}_{(a)}\circ \mathrm{T}_{(a)}=\mathrm{id}$. For this it is enough to show that 
\begin{equation*}
\pi_{(a)}\circ \mathrm{T}_{(a)}g=g\circ\pi_D.
\end{equation*} 
From (\ref{eqn:lift}) we have
\begin{equation*}
\pi_{(a)}\circ \mathrm{T}_{(a)}g(u)=\pi_{(a)}\rho_{(a)}g(\pi_D(u))=g(\pi_D(u)).
\end{equation*}
If $k(a)=0$, then by (\ref{eqn:lift}), $\mathrm{T}_{(a)}f(u)=\rho_{(a)}f(\pi_D(u))$, so 
\begin{equation*}
\mathrm{T}_{(a)}\mathrm{H}_{(a)}g(u)=\rho_{(a)}\mathrm{H}_{(a)}g(\pi_D(u))=\rho_{(a)}\pi_{(a)}g(u)=g(u),
\end{equation*} as $\rho_{(a)}\pi_{(a)}=\mathrm{id}$.

If now $(b)\subset (a)$ and $k(b)<k(a)$, then $b=za$ for some $z\in \mathfrak{m}$. Since $\mathrm{T}_{(b)}g\in (b)$, we see that $\pi_{(a)}\mathrm{T}_{(b)}g=0$. Therefore $\mathrm{H}_{(a)}\circ \mathrm{T}_{(b)}=0$.
\end{proof}

\begin{remark}
Assume that $A$ is a local Frobenius algebra and $B$ is an arbitrary $A$-module. Then Lemma \ref{L:smoothing}, for $k=0$, and Theorem \ref{T:lifting}, for $a=1$, hold. Therefore, we see that if $A/\mathfrak{m}=\mathbb{R}$, then there are $A$-differentiable functions on $B$ which are not smooth and not analytic. If $A$ is not local, then it is a product of local Frobenius algebras, so also the general case can be deduced. Thus for such algebras there exist $A$-differentiable functions which are not $A$-analytic.
\end{remark}

\begin{definition}\label{D:structure}
Assume that $B$ is a finitely generated $A$-module. Let $M\subset \mathbb{N}^r$ be a set of multi-indices, $e_1,\dotsc,e_r\in\mathfrak{m}$ be such as in Lemma \ref{L:generate}. Assume that for any $\alpha\in M$, $B$ and $(e^{\alpha})$ have the lifting property with decomposition $B=C_k\oplus D_k$, $k=\abs{\alpha}$. 
Let $U$ be a convex, open, bounded subset of $B$.
Define
\begin{equation*}
\mathrm{T}\colon \bigoplus_{k=0}^l\bigoplus_{\alpha\in M,\abs{\alpha}=k} C^{p+k(e^{\alpha})}_{A/\mathfrak{m}}(\overline{\pi_{D_k}(U)}, (e^{\alpha})/\mathfrak{m}(e^{\alpha}))\to C^p_A(\overline{U},A)
\end{equation*}
by the formula
\begin{equation*}
\mathrm{T}((f_{\alpha})_{\alpha\in M})=\sum_{\alpha\in M} \mathrm{T}_{(e^{\alpha})}f_{\alpha}.
\end{equation*}
Define
\begin{equation*}
\mathrm{H}\colon C^p_A(\overline{U},A)\to  \bigoplus_{k=0}^l\bigoplus_{\alpha\in M,\abs{\alpha}=k} C^{p+k(e^{\alpha})}_{A/\mathfrak{m}}(\overline{\pi_{D_k(U)}}, (e^{\alpha})/\mathfrak{m}(e^{\alpha}))
\end{equation*}
by the formula
\begin{equation*}
\mathrm{H}f=(f_{\alpha})_{\alpha\in M},
\end{equation*}
where
\begin{equation*}
\begin{aligned}
f_{(0,\dotsc,0)}&=\mathrm{H}_{(e^0)}(f),\\
\sum_{\mathclap{\alpha\in M, \abs{\alpha}=k}}f_{\alpha}[e^{\alpha}]&=\mathrm{H}_{\mathfrak{m}^k}(f-\sum_{\mathclap{\beta\in M, \abs{\beta} < k}}\mathrm{T}_{(e^{\beta})}f_{\beta}),
\end{aligned}
\end{equation*}
for $k=1,\dotsc,l$.
\end{definition}

\begin{remark}
If $f\in C^p_A(\overline{U},\mathfrak{m}^k)$, then from Theorem \ref{T:lifting} it follows that
\begin{equation}\label{eqn:inv1}
\sum_{\mathclap{\alpha\in M, \abs{\alpha}=k}}\mathrm{H}_{\mathfrak{m}^k}\mathrm{T}_{(e^{\alpha})}f_{\alpha}[e^{\alpha}]=\sum_{\mathclap{\alpha\in M, \abs{\alpha}=k}}\mathrm{H}_{(e^{\alpha})}\mathrm{T}_{(e^{\alpha})}f_{\alpha}[e^{\alpha}]=\sum_{\mathclap{\alpha\in M, \abs{\alpha}=k}}f_{\alpha}[e^{\alpha}].
\end{equation}
For $f\in C^p_A(\overline{U},\mathfrak{m}^l)$ we have
\begin{equation}\label{eqn:inv2}
\sum_{\mathclap{\alpha\in M, \abs{\alpha}=l}}\mathrm{T}_{(e^{\alpha})}\mathrm{H}_{(e^{\alpha})}f=f.
\end{equation}
\end{remark}

\begin{theorem}\label{T:structure}
Operators $\mathrm{H}$ and $\mathrm{T}$ are isomorphisms of Banach spaces and are mutual reciprocals.
\end{theorem}
\begin{proof}
Observe that by Lemma \ref{L:propertyL} $B$ and $\mathfrak{m}^k$ have the lifting property for any $k=0,1,\dotsc,l$.
Let us check that $\mathrm{H}$ is well defined. For this we have to show that for $f\in C^p_A(\overline{U},A)$ and all $k=0,1,\dotsc,l$
\begin{equation}\label{eqn:max}
f-\sum_{\mathclap{\beta\in M, \abs{\beta} < k}}\mathrm{T}_{(e^{\beta})}f_{\beta}\in C^p_A(\overline{U},\mathfrak{m}^k).
\end{equation}
For $k=0$ this is the assumption. We shall proceed inductively. Assume that condition (\ref{eqn:max}) holds for $k<l$. Then functions $(f_{\alpha})_{\alpha\in M,\abs{\alpha}=k}$ are well defined by the formula
\begin{equation*}
\sum_{\mathclap{\alpha\in M, \abs{\alpha}=k}}f_{\alpha}[e^{\alpha}]=\mathrm{H}_{\mathfrak{m}^k}(f-\sum_{\mathclap{\beta\in M, \abs{\beta} < k}}\mathrm{T}_{(e^{\beta})}f_{\beta}).
\end{equation*}
Moreover, by (\ref{eqn:inv1}), we have
\begin{equation*}
\mathrm{H}_{\mathfrak{m}^k}(f-\sum_{\mathclap{\beta\in M, \abs{\beta} < k+1}}\mathrm{T}_{(e^{\beta})}f_{\beta})=\sum_{\mathclap{\alpha\in M, \abs{\alpha}=k}}f_{\alpha}[e^{\alpha}]-\sum_{\mathclap{\beta\in M, \abs{\beta} = k}}\mathrm{H}_{(e^{\beta})}(\mathrm{T}_{(e^{\beta})}f_{\beta})[e^{\beta}]=0.
\end{equation*}
This completes the induction.
We shall now check that $\mathrm{T}$ and $\mathrm{H}$ are mutual reciprocals. Let $f\in C^p_A(\overline{U},A)$ and $\mathrm{H}f=(f_{\alpha})_{\alpha\in M}$. Then, by (\ref{eqn:inv2}),
\begin{equation*}
\sum_{\mathclap{\beta\in M, \abs{\beta} = l}}\mathrm{T}_{(e^{\beta})}f_{\beta}=\sum_{\mathclap{\beta\in M, \abs{\beta} = l}}\mathrm{T}_{(e^{\beta})}\mathrm{H}_{(e^{\beta})}(f-\sum_{\mathclap{\alpha\in M, \abs{\alpha} < l}}\mathrm{T}_{(e^{\alpha})}f_{\alpha})
=f-\sum_{\mathclap{\alpha\in M, \abs{\alpha} < l}}\mathrm{T}_{(e^{\alpha})}f_{\alpha}.
\end{equation*}
Therefore
\begin{equation*}
f=\sum_{\mathclap{\alpha\in M}}\mathrm{T}_{(e^{\alpha})}f_{\alpha}.
\end{equation*}
That is $\mathrm{id}=\mathrm{T}\circ \mathrm{H}$.
Choose now some
\begin{equation*}
f\in \bigoplus_{k=0}^l\bigoplus_{\alpha\in M,\abs{\alpha}=k} C^{p+k(e^{\alpha})}_{A/\mathfrak{m}}(\overline{\pi_{D_k}(U)}, (e^{\alpha})/\mathfrak{m}(e^{\alpha})). 
\end{equation*}
Then $\mathrm{T}f=\sum_{\alpha\in M}T_{(e^{\alpha})}f_{\alpha}$. By Theorem \ref{T:lifting} we have $\mathrm{H}_{(e^0)}\mathrm{T}f=f_{0}$. Suppose now that for  $k<l$ and all $\alpha\in M$ such that $\abs{\alpha}\leq k$ we have $(\mathrm{HT}f)_{\alpha}=f_{\alpha}$. Then, by the definition of $\mathrm{T}$,
\begin{equation*}
\mathrm{T}f-\sum_{\mathclap{\alpha\in M, \abs{\alpha}\leq k}}\mathrm{T}_{(e^{\alpha})}(\mathrm{HT}f)_{\alpha}=\sum_{\mathclap{\beta\in M, \abs{\beta}>k}}\mathrm{T}_{(e^{\beta})}f_{\beta}.
\end{equation*}
By the definition of $\mathrm{H}$ we have
\begin{equation*}
\begin{aligned}
&\sum_{\mathclap{\gamma\in M, \abs{\gamma}=k+1}}(\mathrm{HT}f)_{\gamma}[e^{\gamma}]=\mathrm{H}_{\mathfrak{m}^{k+1}}(\mathrm{T}f-\sum_{\mathclap{\alpha\in M,\abs{\alpha}\leq k}}\mathrm{T}_{(e^{\alpha})}(\mathrm{HT}f)_{\alpha})=\\
&=\sum_{\mathclap{\beta\in M, \abs{\beta}>k}}\mathrm{H}_{\mathfrak{m}^{k+1}}(\mathrm{T}_{(e^{\beta})}f_{\beta})=\sum_{\mathclap{\beta\in M, \abs{\beta}=k+1}}\mathrm{H}_{\mathfrak{m}^{k+1}}(\mathrm{T}_{(e^{\beta})}f_{\beta})=\sum_{\mathclap{\beta\in M, \abs{\beta}=k+1}}f_{\beta}[e^{\beta}].
\end{aligned}
\end{equation*}
From induction it follows that $\mathrm{H}\circ \mathrm{T}=\mathrm{id}$.
Continuity of $\mathrm{T}$ and $\mathrm{H}$ follows from Lemma \ref{L:smoothing} and Theorem \ref{T:lifting}.
\end{proof}

\begin{remark}
Since $A/\mathfrak{m}$ is a finite-dimensional extension of $\mathbb{F}$, there are two cases to consider. Either $A/\mathfrak{m}$ is equal to $\mathbb{R}$ or to $\mathbb{C}$. We want to stress the difference which occurs.

In the complex field case, differentiability of a function in  
\begin{equation*}
C^{p+k(a)}_{A/\mathfrak{m}}(\overline{\pi_D(U)},(a)/\mathfrak{m}(a))
\end{equation*}
in open set $\pi_D(U)$ is automatic for any $p\geq 0$. Indeed, the condition of $A/\mathfrak{m}$-continuity means that a function is complex differentiable. However, the continuity of derivatives on the closure of $U$ is not automatic. 

In the real field case, which appears to be much more interesting, we see that certain components of a function in $C^p_A(\overline{U},(a))$ are necessarily of higher differentiability, which is an unexpected phenomenon.
\end{remark}

\section{General algebra, free module}

We shall now consider the case when $A$ is an arbitrary finite dimensional, commutative algebra and the module $B$ is free. We may, without loss of generality, assume that $A$ is local. The maximal ideal in $A$ is denoted by $\mathfrak{m}$.

Let $b_1,\dotsc,b_n$ be an $A$-basis of $B$. Define $\rho_B\colon B/\mathfrak{m}B\to B$ by the formula
\begin{equation*}
\rho_D(\sum_{i=1}^n[a_i]\pi_B(b_i))=\sum_{i=1}^n [a_i]b_i,
\end{equation*}
for any $a_i\in A$, $i=1,\dotsc,n$.

\begin{lemma}\label{L:freemod}
Assume that $B$ is a finitely generated, free $A$-module. Then for any $a\in A$, $B$ and $(a)$ have the lifting property with the decomposition $B=C\oplus B$, $C=\{0\}$. 
\end{lemma}
\begin{proof}
Let $a\in A$. We have to show that any symmetric $A/\mathfrak{m}$-multilinear map $f\colon B\times\dotsm\times B\to (a)/\mathfrak{m}(a)$ admits a lifting to an $A$-multilinear symmetric map $h\colon B\times\dotsm\times B\to (a)$ such that $\pi_{(a)}\circ h=f$.
Let $f\colon B\times\dotsm\times B\to (a)/\mathfrak{m}(a)$ be symmetric and $A/\mathfrak{m}$-multilinear. Define 
\begin{equation*}
\begin{aligned}
\mathrm{G}^m_{(a)}f&(\sum_{i_1=1}^na_{i_11}b_{i_1},\dotsc,\sum_{i_m=1}^na_{i_mm}b_{i_m})=\\
&=\sum_{i_1=1,\dotsc,i_m=1}^n a_{i_11}\dotsm a_{i_mm} \rho_{(a)}f(b_{i_1},\dotsc,b_{i_m}).
\end{aligned}
\end{equation*}
Then conditions (\ref{i:comp})-(\ref{i:comp1}) and (\ref{i:comp3}) of Definition \ref{D:L} are clearly satisfied. Let us check that the condition (\ref{i:comp2}) holds true.
Let $f\colon B\to (a)/\mathfrak{m}(a)$ be $A/\mathfrak{m}$-linear. Recall that $\pi_D=\pi_B$ is a quotient map onto $B/\mathfrak{m}B$ and $\rho_D=\rho_B$ is $A/\mathfrak{m}$-linear map from $B/\mathfrak{m}B$ to $B$. We have to check that
\begin{equation*}
\mathrm{G}^1_{(a)}f(\sum_{i=1}^n [a_i]b_i)=\rho_{(a)}f(\sum_{i=1}^n [a_i]b_i).
\end{equation*}
This is again clear from the definition of $\mathrm{G}^1_{(a)}$.
\end{proof}

\begin{theorem}\label{T:structurefree}
Assume that $B$ is a finitely generated, free module over a local algebra $A$. Assume that $U\subset B$ is open, convex and bounded. Then any $A$-continuous function in $f\in C^p_A(\overline{U},A)$ may be written in the form
\begin{equation*}
f=\sum_{\alpha\in M}\left(\rho_{(e^{\alpha})}f_{\alpha}(\pi_B(u))+\sum_{j=1}^{k(e^{\alpha})}\frac1{j!}\mathrm{G}^j_{(e^{\alpha})}\mathrm{D}^jf_{\alpha}(\pi_B(u))((u-\rho_B\pi_B(u))^j)\right).
\end{equation*}
for some functions $(f_{\alpha})_{\alpha\in M}\in  \bigoplus_{\alpha\in M} C^{p+k(e^{\alpha})}_{A/\mathfrak{m}}(\overline{\pi_B(U)}, (e^{\alpha})/\mathfrak{m}(e^{\alpha}))$
Conversely, any such function belongs to $C^p_A(\overline{U},A)$. This assignment is an isomorphism of Banach spaces.
\end{theorem}
\begin{proof}
Follows directly from Lemma \ref{L:freemod} and Theorem \ref{T:structure}.
\end{proof}

This theorem answers the question raised by Waterhouse in \cite{W2} about rules satisfied by $A$-differentiable functions on algebras. The case of an arbitrary finitely generated module remains open, but it can be expected that the structure of $A$-differentiable functions should rely on the lifting properties of $A$-module homomorphisms. Another possible approach is to follow the observation from Lemma \ref{L:etadifferentiability} that if $f\in C^p_A(\overline{U},A)$, and $\eta\colon A^n\to B$ is a surjective $A$-linear map, then $f\circ \eta$ is in $C^p_A(\overline{\eta^{-1}(U)},A)$. Thus Theorem \ref{T:structurefree} gives us a description of $f\circ \eta$. 

\section{Algebra generated by one element}\label{S:onegenerator}

We shall now study the case when $A$ is an algebra generated by one element. This case is particularly simple and will lead us to interesting applications. Let us first describe how such algebras look and how modules over such algebras look.
We refer to \cite{JAC2, JAC} to background on algebras and their modules.

Any finite dimensional, commutative algebra over $\mathbb{F}$ that is generated by one element is isomorphic to $\mathbb{F}[x]/(P(x))$, for some polynomial $P$.

If $\mathbb{F}=\mathbb{R}$, then $P$ factors into a product of polynomials of the form\footnote{We shall use this notation also for $\lambda\in \mathbb{C}$.}
\begin{equation}\label{eqn:firstf}
Q_{\lambda,k}(x)=(x-\lambda)^k,
\end{equation}
for some $\lambda\in\mathbb{R}$ or of the form
\begin{equation}\label{eqn:secondf}
R_{\alpha,\beta,k}(x)=(x^2-2\alpha x+\beta)^k,
\end{equation}
for some $\alpha,\beta\in\mathbb{R}$ such that $\alpha^2+\beta>0$. That is
\begin{equation*}
P(x)=\prod_{i=1}^{m} Q_{\lambda_i,l_i}(x) \prod_{i=1}^{n}R_{\alpha_i,\beta_i,k_i}(x).
\end{equation*}

Sending $x\mapsto z$ we have $\mathbb{R}[x]/(R_{\alpha_i,\beta_i,k_i}(x))\cong \mathbb{C}[z]/(Q_{\gamma_i,k_i}(z))$, where $\gamma_i=\alpha_i\pm \sqrt{\alpha_i^2+\beta_i} i$. 

By the Chinese remainder theorem we obtain a decomposition of $\mathbb{R}[x]/(P(x))$ into a direct sum of local algebras
\begin{equation}\label{eqn:real}
\mathbb{R}[x]/(P(x))\cong \bigoplus_{i=0}^{m}\mathbb{R}[x]/(Q_{\lambda_i,l_i}(x))\oplus \bigoplus_{i=0}^n \mathbb{C}[z]/(Q_{\gamma_i,k_i}(z)).
\end{equation}

When $\mathbb{F}=\mathbb{C}$, then the polynomial $P$ factors into a product of polynomials of the form (\ref{eqn:firstf}). Thus
\begin{equation}\label{eqn:complex}
\mathbb{C}[x]/(P(x))\cong \bigoplus_{i=1}^n \mathbb{C}[z]/(Q_{\gamma_i,k_i}(z)).
\end{equation}

Proposition \ref{P:spliting} tells us that we may restrict ourselves to the case 
\begin{equation*}
A=\mathbb{F}[x]/(Q_{\lambda,l}(x)).
\end{equation*}
Let $e=Q_{\lambda,1}(x)$. Then $e$ generates the maximal ideal $\mathfrak{m}=(Q_{\lambda,1}(x))$. Note that $\mathfrak{m}^{l-1}\neq\{0\}$ and $l-1$ is maximal among all such natural numbers. Moreover powers of $e$, including $1=e^0$, are an $A/\mathfrak{m}$-basis of $A$.

From the structure theorem for finitely generated modules over principal ideal domains (see \cite{JAC}) or from the Jordan canonical form of a matrix, any finitely generated module over such an algebra is of the form
\begin{equation*}
B=\bigoplus_{i=1}^{j_1} \mathbb{F}[x]/(Q_{\lambda,1}(x))\oplus \bigoplus_{i=1}^{j_2} \mathbb{F}[x]/(Q_{\lambda,2}(x))\oplus\dots\oplus \bigoplus_{i=1}^{j_l} \mathbb{F}[x]/(Q_{\lambda,l}(x)),
\end{equation*}
for some natural numbers $j_1,\dotsc,j_l\geq 0$; if $j_t=0$ for some $t=1,\dotsc,l$, then the corresponding summand is omitted.

Let $e_t$ denote the unit in the summand $\mathbb{F}[x]/(Q_{\lambda,t}(x))$. Let for $k=0,1,\dotsc,l-1$
\begin{equation*}
C_k=\bigoplus_{i=1}^{j_1} \mathbb{F}[x]/(Q_{\lambda,1}(x))\oplus\dots\oplus \bigoplus_{i=1}^{j_{l-k-1}} \mathbb{F}[x]/(Q_{\lambda,l+k-1}(x)),
\end{equation*}
and
\begin{equation*}
D_k=\bigoplus_{i=1}^{j_{l-k}} \mathbb{F}[x]/(Q_{\lambda,l-k}(x))\oplus\dots\oplus \bigoplus_{i=1}^{j_l} \mathbb{F}[x]/(Q_{\lambda,l}(x)).
\end{equation*}
Let  
\begin{equation*}
\rho_{(e^k)}\colon (e^k)/\mathfrak{m}(e^k)\to (e^k)
\end{equation*}
be defined by the formula $\rho_{(e^k)}(a[e^k])=ae^k$ for $a\in A/\mathfrak{m}$. Then $\pi_{(e^k)}\circ \rho_{(e^k)}=\mathrm{id}$.
Let 
\begin{equation*}
\rho_{D_k}\colon D_k/\mathfrak{m}D_k\to B
\end{equation*}
be defined by the formula $\rho_{D_k}(a[e_t])=ae_t$ for $a\in A/\mathfrak{m}$ and $t\geq l-k$. Then $\pi_{D_k}\circ \rho_{D_k}=\mathrm{id}$.

\begin{lemma}\label{L:onegeneratorL}
Let $k=0,\dotsc,l-1$. Then the module $B$ and ideal $(e^k)\subset A$ have the lifting property with the decomposition $B=C_k\oplus D_k$.
\end{lemma}
\begin{proof}
Assume that for an $A/\mathfrak{m}$-linear map $f\colon B\to  (e^k)/\mathfrak{m}(e^k)$ there exists an $A$-linear map $h\colon B\to (e^k)$, such that $\pi_{(e^k)}\circ h=f$. Then 
\begin{equation*}
0=h(Q_{\lambda,t}(x)e_t)=Q_{\lambda,t}(x)h(e_t).
\end{equation*}
Since $(Q_{\lambda,t}(x))=\mathfrak{m}^t$, we see that $h(e_t)\in Ann(\mathfrak{m}^t)=\mathfrak{m}^{l-t}$. Thus if $l-t\geq k+1$, that is $t\leq l-k-1$, then $\pi_{(e^k)}\circ h(e_t)=f(e_t)=0$. As $e_t$ generate $\mathbb{F}[x]/(Q_{\lambda,t}(x))$, $f$ is zero on whole $C_k$.

Assume now that $f\colon B\times\dotsm\times B \to (e^k)/\mathfrak{m}(e^k)$ is a multilinear map such that for any $b_1,\dotsc,b_{j-1}\in B$ 
\begin{equation*}
f(b_1,\dotsc,b_{m-1},\cdot,b_{m+1},\dotsc,b_{j-1})|_{C_k}=0.
\end{equation*}
Set 
\begin{equation}\label{eqn:rho}
h(e_{t_1},\dotsc,e_{t_j})=\rho_{(e^k)}f(e_{t_1},\dotsc,e_{t_j})
\end{equation}
for all $t=1,\dotsc,l$ and all $i=1,\dotsc,j_t$. Extend this definition by $A$-linearity
\begin{equation}\label{eqn:extension}
h(Q_{\lambda,s_1}(x)e_{t_1},\dotsc,Q_{\lambda,s_j}(x)e_{t_j})=Q_{\lambda,s_1}(x)\dotsm Q_{\lambda,s_j}(x)\rho_{(e^k)}f(e_{t_1},\dotsc,e_{t_j}),
\end{equation}
for all $s_r<t_r$, $r=1,\dotsc,j$. Since $Q_{\lambda,t}(x)e_t=0$, we must have 
\begin{equation*}
Q_{\lambda,t}(x)\rho_{(e^k)}f(e_{t_1},\dotsc,e_t,\dotsc,e_{t_j})=0.
\end{equation*} 
By assumption $f(e_{t_1},\dotsc,e_t,\dotsc,e_{t_j})=0$ if $t\leq l-k-1$. Assume that $t\geq l-k$. Then $\rho_{(e^k)}f(e_{t_1},\dotsc,e_t,\dotsc,e_{t_j})\in\mathfrak{m}^k$ and $Q_{\lambda,t}(x)\in\mathfrak{m}^t$, so 
\begin{equation*}
Q_{\lambda,t}(x)\rho_{(e^k)}f(e_{t_1},\dotsc,e_t,\dotsc,e_{t_j})\in \mathfrak{m}^{t+k}\subset\mathfrak{m}^l=\{0\}.
\end{equation*}
This shows that $h$ is well defined. 

Furthermore, we have
\begin{equation*}
(\pi_{(e^k)}\circ h) (e_{t_1},\dotsc,e_{t_j})=(\pi_{(e^k)}\circ \rho_{(e^k)})f(e_{t_1},\dotsc,e_{t_j})=f(e_{t_1},\dotsc,e_{t_j}).
\end{equation*}
Since both sides are $A$-linear and $e_t$ generate $B$, $\pi_{(e^k)}\circ h=f$. 

We define $\mathrm{G}^j_{(e^k)}f=h$. Then $\mathrm{G}^j_{(e^k)}$ are linear. Choose a function $f\colon B\to (e^k)/\mathfrak{m}(e^k)$ which vanishes on $C_k$ and $v\in B$. 
It is necessary now to show that
\begin{equation*}
\mathrm{G}^j_{(e^k)}f(\rho_{D_k}\pi_{D_k}(v))=\rho_{(e^k)}f(\rho_{D_k}\pi_{D_k}(v)).
\end{equation*}
This follows immediately from (\ref{eqn:rho}).
The condition (\ref{i:comp3}) from Definition \ref{D:L} follows readily - both sides there are linear, so it is enough to check it on a basis. Thus it follows from (\ref{eqn:extension}).
\end{proof}

\begin{theorem}\label{T:structure-one}
Assume that $B$ is a finitely generated module over $A=\mathbb{F}[x]/(Q_{\lambda,l})$. Assume that $U\subset B$ is an open, convex and bounded set. Then any $A$-continuous function in $f\in C^p_A(\overline{U},A)$ may be written in the form
\begin{equation*}
f=\sum_{k=0}^{l-1}\left(\rho_{(e^k)}f_k(\pi_{D_k}(u))+\sum_{j=1}^{l-1-k}\frac1{j!}\mathrm{G}^j_{(e^k)}\mathrm{D}^jf_k(\pi_{D_k}(u))((u-\rho_{D_k}\pi_{D_k}(u))^j)\right).
\end{equation*}
for some functions 
\begin{equation*}
(f_k)_{k=0,1\dotsc,l-1}\in  \bigoplus_{k=0}^{l-1} C^{p+l-1-k}_{A/\mathfrak{m}}(\overline{\pi_{D_k}(U)}, (e^k)/\mathfrak{m}(e^k)).
\end{equation*}
Conversely, any such function belongs to $C^p_A(\overline{U},A)$. This assignment is an isomorphism of Banach spaces.
\end{theorem}
\begin{proof}
Follows readily from Lemma \ref{L:onegeneratorL} and Theorem \ref{T:structure}.
\end{proof}

\begin{lemma}\label{L:froboundary}
Let $A=\mathbb{F}[x]/(Q_{\lambda,l}(x))$. Let $\phi\colon A\to \mathbb{F}$ be defined by the formula
\begin{equation*}
\phi(\sum_{i=0}^{l-1}a_ie^i)=a_{l-1},
\end{equation*}
where $a_i\in \mathbb{F}$.
Assume that $\overline{\rho_{D_i}\pi_{D_i}(U)}\subset \overline{U}$ for all $i=0,\dotsc,l-1$. Then for any $k=0,\dotsc,l-1$ and any 
\begin{equation*}
f\in C^{p+l-1-k}_{A/\mathfrak{m}}(\overline{\pi_{D_k}(U)}, (e^k)/\mathfrak{m}(e^k))
\end{equation*}
we have
\begin{equation*}
\phi(e^{l-1-k}\mathrm{T}_{(e^k)}f)|_{\overline{\rho_{D_k}\pi_{D_k}(U)}}=f\circ \pi_{D_k}
\end{equation*}
and if $j\neq k$
\begin{equation*}
\phi(e^{l-1-j}\mathrm{T}_{(e^k)}f)|_{\overline{\rho_{D_j}\pi_{D_j}(U)}}=0.
\end{equation*}
If functions 
\begin{equation*}
f_k\in C^{p+l-1-k}_{A/\mathfrak{m}}(\overline{\pi_{D_k}(U)}, (e^k)/\mathfrak{m}(e^k)),
\end{equation*}
$k=0,\dotsc,l-2$, have derivatives equal to zero and 
\begin{equation*}
f_{l-1}\in C^p_{A/\mathfrak{m}}(\overline{\pi_{D_{l-1}}(U)}, (e^{l-1})/\mathfrak{m}(e^{l-1}))
\end{equation*}
is equal to zero, then 
\begin{equation*}
\phi(\mathrm{T}f)=0.
\end{equation*}
\end{lemma}
\begin{proof}
By (\ref{eqn:lift}) it follows that $\mathrm{T}_{(e^k)}f|_{\overline{\rho_{D_k}\pi_{D_k}(U)}}=\rho_{(e^k)}f\circ \pi_{D_k}$. Then
\begin{equation*}
\phi(e_{l-1-k}\mathrm{T}_{(e^k)}f)|_{\overline{\rho_{D_k}\pi_{D_k}(U)}}=\phi(e_{l-1-k}f\circ \pi_{D_k}e_k)|_{\overline{\rho_{D_k}\pi_{D_k}(U)}}=f\circ \pi_{D_k}.
\end{equation*}
and if $j\neq k$
\begin{equation*}
\phi(e_{l-1-j}\mathrm{T}_{(e^k)}f)|_{\overline{\rho_{D_k}\pi_{D_k}(U)}}=\phi(e_{l-1-j}f\circ \pi_{D_k}e_k)|_{\overline{\rho_{D_k}\pi_{D_k}(U)}}=0.
\end{equation*}
Let $f_k\in C^{p-k+l-1}_{A/\mathfrak{m}}(\overline{\rho_{D_k}\pi_{D_k}(U)}, (e^k)/\mathfrak{m}(e^k))$, $k=0,\dotsc,l-2$ have derivatives equal to zero and let $f_{l-1}\in C^p_{A/\mathfrak{m}}(\overline{\rho_{D_{l-1}}\pi_{D_{l-1}}(U)}, (e^{l-1})/\mathfrak{m}(e^{l-1}))$ be equal to zero. Then by (\ref{eqn:lift}) it follows that $\mathrm{T}_{(e^k)}f_k=\rho_{(e^k)}f_k\circ \pi_{D_k}$ for $k=0,\dotsc,l-2$ and $\mathrm{T}_{l-1}f_{l-1}=0$.
Since $\phi(e^k)=0$ for $k=0,\dotsc,l-2$ we have
\begin{equation*}
\phi(\mathrm{T}f)=\sum_{k=0}^{l-2}\phi(\rho_{(e^k)}f_k\circ \pi_{D_k})=0.
\end{equation*} 
\end{proof}

\section{The equation $grad(w)=\tens{M}grad(v)$}

This section is devoted to application of the theory to the equation $grad(w)=\tens{M}grad(v)$. In particular, we study the case, when the matrix $\tens{M}$ has at least two Jordan blocks corresponding to the same eigenvalue, which was not covered by Waterhouse in \cite{W1}.

Our equation can be reformulated in an equivalent form
\begin{equation*}
\mathrm{D}w(\cdot)(x)=\mathrm{D}v(\cdot)(\tens{M}^{\intercal}x),
\end{equation*}
for all $x\in \mathbb{F}^n$.

Observe that an algebra $A$ generated by a matrix $M$ is isomorphic to $\mathbb{F}[x]/(P(x))$, where $P$ is the minimal polynomial of $M$. By Example \ref{E:Frobenius} $A$ is a Frobenius algebra.

We prove an analogue of Theorem 5.1 from \cite{W1}. 

\begin{theorem}\label{T:pairs}
Let $t\geq 0$ be a natural number. Let $U\subset \mathbb{F}^n$ be an open, simply connected set.
Let $A$ be an algebra generated by the matrix $\tens{M}^{\intercal}\in M_{n\times n}(\mathbb{F})$ and let $B=\mathbb{F}^n$ with the natural structure of an $A$-module. Let $v, w\colon U \to \mathbb{F}$. The following conditions are equivalent
\begin{enumerate}[\upshape (i)]
\item $v$, $w$ are $C^{2+t}(U)$ functions satisfying $\mathrm{D}w=\mathrm{D}v\circ \tens{M}^{\intercal}$,
\item $v=\phi(f)$ is a component function of an $A$-differentiable function $f\colon U \to A$ of class $C^{2+t}(U,A)$, and $w=\phi(\tens{M}^{\intercal}f)+c$, where $c$ is a constant.
\end{enumerate}
Assume additionally that $U$ is short-path connected\footnote{See Definition \ref{D:short-path}.}. Then the following conditions are equivalent
\begin{enumerate}[\upshape (i)]
\item $v$, $w$ are $C^{2+t}(\overline{U})$ functions satisfying $\mathrm{D}w=\mathrm{D}v\circ \tens{M}^{\intercal}$,
\item $v=\phi(f)$ is a component function of an $A$-differentiable function $f\colon U \to A$ of class $C^{2+t}(\overline{U},A)$, and $w=\phi(\tens{M}^{\intercal}f)+c$, where $c$ is a constant.
\end{enumerate} 
\end{theorem}
\begin{proof}
Let us first prove the first equivalence. Assume that $f$ is $A$-differentiable and of class $C^{2+t}(U,A)$ and $v=\phi(f)$ and $w=\phi(\tens{M}^{\intercal}f)+c$ . Then, thanks to $A$-linearity of the derivative of $f$, we have: 
\begin{equation*}
\mathrm{D}w(b)(x)=\phi(\mathrm{D}\tens{M}^{\intercal}f(b)(x))=\phi(\mathrm{D}f(b)(\tens{M}^{\intercal}x))=\mathrm{D}v(b)(M^{\intercal}x).
\end{equation*} 
Conversely, if $v$ and $w$ satisfy $\mathrm{D}w=\mathrm{D}v\circ \tens{M}^{\intercal}$, then the right-hand side of the equation is a derivative of some $C^{2+t}(U)$ function. Hence it's second derivative must be symmetric. Thus 
\begin{equation*}
\begin{aligned}
\mathrm{D}^2v(b)(\tens{M}^{\intercal}x,y)&=\mathrm{D}(\mathrm{D}v(b)(\tens{M}^{\intercal}x))(y)=\mathrm{D}(\mathrm{D}w(b)(x))(y)=\\
&=\mathrm{D}^2w(b)(x,y)=D^2w(b)(y,x)=D^2v(b)(\tens{M}^{\intercal}y,x).
\end{aligned}
\end{equation*}
As $\tens{M}^{\intercal}$ generates the algebra, this condition holds also for any element $a\in A$. By Theorem \ref{T:components}, there exists an uniquely determined, up to a constant, $A$-differentiable $f\colon U \to A$, of class $C^{2+t}(U,A)$, which satisfies $v=\phi(f)$. Then $w$ and $\phi(\tens{M}^{\intercal}f)$ have the same derivative by the first part of the proof, so they differ by a constant.

The second equivalence follows from Theorem \ref{T:components} in the same way.  
\end{proof}

\begin{theorem}\label{T:analyticity}
Let $U$ be an open, convex and bounded subset of $\mathbb{F}^n$. The following conditions are equivalent:
\begin{enumerate}[\upshape (i)]
\item any functions $v$, $w$ of class $C^2(U)$ which satisfy $\mathrm{D}w=\mathrm{D}v\circ \tens{M}^{\intercal}$, are analytic and are components of $A$-analytic functions.
\item $\tens{M}$ has no real eigenvalues
\end{enumerate}
\end{theorem}
\begin{proof}
If $\tens{M}$ does not have any real eigenvalue, then the algebra $A$ generated by $\tens{M}^{\intercal}$ is a $\mathbb{C}$-algebra. By Theorem \ref{T:C-continuous} every $A$-differentiable function is $A$-analytic. So is a component of such a function. 

Conversely, if $\tens{M}$ has a real eigenvalue, then by Proposition \ref{P:spliting} and \S \ref{S:onegenerator} we may assume that it is the only eigenvalue. Then Theorem \ref{T:structure-one} gives us examples of $A$-differentiable functions whose all components are not analytic.
\end{proof}

\section{Boundary value problem for generalised Laplace equations}\label{S:boundary}

Now we are ready to pose a boundary value problem. We assume that $U$ is an open, convex, bounded set. Our aim is to find appropriate boundary for the equation and the conditions that have to be imposed on the boundary values so that a solution is of prescribed differentiability and is unique up to a constant. Theorem \ref{T:pairs} tells us that solving the equation, up to a constant, is the same as finding an $A$-differentiable function. By Theorem \ref{T:structure-one} we see that to describe an $A$-differentiable function it is sufficient and enough to prescribe some functions of sufficient differentiability on the projections of $U$ and then extend them to the whole of $U$ by the Taylor's formula. 

\begin{example}
Consider the matrix 
\begin{equation*}
\tens{M}=\begin{bmatrix}
\lambda&0&0\\
0&\lambda&0\\
0&1&\lambda
\end{bmatrix}.
\end{equation*}
The minimal polynomial of $\tens{M}^{\intercal}$ is $P(x)=(x-\lambda)^2$. Algebra $A$ generated by $\tens{M}^{\intercal}$ is then equal to $\mathbb{R}[x]/(x-\lambda)^2$. $A$-module $B=\mathbb{R}^3$ has the decomposition 
\begin{equation*}
B=\mathbb{R}[x]/(x-\lambda)\oplus\mathbb{R}[x]/(x-\lambda)^2.
\end{equation*}
$A$ has a basis $1,e$, where $e=x-\lambda$. $B$ has a basis $(e_1,e_2,e_3)$, such that $ee_1=0$, $ee_2=0$, $ee_3=e_2$.  $A$ is a Frobenius algebra, with the functional $\phi(x_11+x_2e)=x_2$, $x_1,x_2\in \mathbb{R}$. Let 
\begin{equation*}
U=(0,1)^3=\{x_1e_1+x_2e_2+x_3e_3\in B\colon 0< x_i<1\}.
\end{equation*}
Let $t\in \mathbb{N}$, $v\in C^{2+t}(\overline{U})$. Consider the equation 
\begin{equation}\label{eqn:dife}
\mathrm{D}^2v(\cdot)(z,\tens{M}^{\intercal}y)=\mathrm{D}^2v(\cdot)(\tens{M}^{\intercal}z,y),
\end{equation}
for all $z,y\in B$.
Equivalently $\mathrm{D}^2v(\cdot)(z,ey)=\mathrm{D}^2v(\cdot)(ez,y)$. This means that 
\begin{equation*}
\mathrm{D}^2v(\cdot)(z_1e_1+z_2e_2+z_3e_3,y_3e_2)=\mathrm{D}^2v(\cdot)(z_3e_2,y_1e_1+y_2e_2+y_3e_3),
\end{equation*}
that is
\begin{equation*}
\frac{\partial^2v}{\partial x_1\partial x_2}=0, \qquad
\frac{\partial^2v}{\partial x_2^2}=0.
\end{equation*}
Theorem \ref{T:components} tells us that any such $v$ is given by $v=\phi(f)$ for some $A$-differentiable $f$ of class $C^2(\overline{U},A)$. Further, by Theorem \ref{T:structure-one}, any such $f$ is uniquely determined by two functions - $f_0$ in $C^3_{\mathbb{R}}(\overline{\pi_{D_0}(U)},\mathbb{R})$, and $f_1$ in $C^2_{\mathbb{R}}(\overline{\pi_{D_1}(U)},\mathbb{R})$, where
\begin{equation*}
\begin{aligned}
&D_0=\mathbb{R}[x]/(x-\lambda)^2\\
&D_1=\mathbb{R}[x]/(x-\lambda)\oplus\mathbb{R}[x]/(x-\lambda)^2,
\end{aligned}
\end{equation*}
and
\begin{equation*}
\begin{aligned}
&\pi_{D_0}\colon B \to D_0/\mathfrak{m}D_0,\quad \pi_{D_0}(x_1e_1+x_2e_2+x_3e_3)=x_2[e_2],\\
&\pi_{D_1}\colon B \to D_1/\mathfrak{m}D_1,\quad \pi_{D_1}(x_1e_1+x_2e_2+x_3e_3)=x_1[e_1]+x_2[e_2].
\end{aligned}
\end{equation*}
The extension is given by
\begin{equation*}
\begin{aligned}
f(u)&=\rho_{(e^0)}f_0(\pi_{D_0}(u))+\mathrm{G}^1_{(e^0)}(\mathrm{D}f_0(\pi_{D_0}(u))(u-\rho_{D_0}\pi_{D_0}u)+\rho_{(e^1)}f_1(\pi_{D_1}(u))=\\
&=f_0(u_2[e_2])1+u_3\frac{\partial f_0}{\partial x_2}(u_2[e_2])e+f_1(u_1[e_1]+u_2[e_2])e.
\end{aligned}
\end{equation*}
Thus 
\begin{equation*}
v(u)=\phi(f(u))=u_3\frac{\partial f_0}{\partial x_2}(u_2[e_2])+f_1(u_1[e_1]+u_2[e_2]).
\end{equation*}
Therefore any solution to the equation $grad(w)=\tens{M}grad(v)$ is of the form
\begin{equation*}
\begin{aligned}
&v(u)=\phi(f(u))=u_3\frac{\partial f_0}{\partial x_2}(u_2[e_2])+f_1(u_1[e_1]+u_2[e_2]),\\
&w(u)=\phi(\tens{M}^{\intercal}f(u))+c=f_0(u_2[e_2])+\lambda u_3\frac{\partial f_0}{\partial x_2}(u_2[e_2])+\lambda f_1(u_1[e_1]+u_2[e_2])+c.
\end{aligned}
\end{equation*}
We see  that there is a unique solution $v$ of generalised Laplace equations (\ref{eqn:dife}) such that it has fixed values on $\overline{\rho_{D_1}\pi_{D_1}(U)}$ and such that $w-\lambda v$ has fixed, up to a constant, values on $\overline{\rho_{D_0}\pi_{D_0}(U)}$.
\end{example}

The last observation in the above example may be generalised to a theorem that links the components of $A$-differentiable functions with the solutions of boundary value problem for generalised Laplace equations.

For clarity, we shall only treat the case when $\tens{M}$ has one eigenvalue. The general case may also be inferred from this particular case, using Proposition \ref{P:spliting}. However, we shall not state the corresponding theorem, due to the complexity of the notation. 

Assume first that the matrix $\tens{M}\in M_{n\times n}(\mathbb{R})$ has one eigenvalue $\lambda\in\mathbb{F}$, $\mathbb{F}\in\{\mathbb{R},\mathbb{C}\}$. Then the algebra $A$ generated by $\tens{M}^{\intercal}$ is a local algebra isomorphic to $\mathbb{F}[x]/(x-\lambda)^l$. Define $\phi_{\mathbb{R}}\colon A\to \mathbb{R}$ by the formula
\begin{equation*}
\phi(\sum_{i=0}^{l-1}a_i[(x-\lambda)^i])=a_{l-1}.
\end{equation*}
and $\phi_{\mathbb{C}}\colon A\to \mathbb{R}$ by the formula
\begin{equation*}
\phi(\sum_{i=0}^{l-1}a_i[(x-\lambda)^i])=Im a_{l-1}.
\end{equation*}
Then $\phi_{\mathbb{F}}$ makes $A$ a Frobenius algebra.
Define $\mu_{\mathbb{F}}\colon\mathbb{F}\to\mathbb{R}$ by $\mu_{\mathbb{R}}=\mathrm{id}$ and $\mu_{\mathbb{C}}=Im$.
Recall that $(x-\lambda)^{k}=e^k$, for all $k=0,\dotsc,l-1$.

\begin{theorem}\label{T:boundary}
Assume that $U\subset\mathbb{R}^n$ is a convex, open and bounded set. Assume that $\overline{\rho_{D_i}\pi_{D_i}(U)}\subset \overline{U}$ for all $i=0,\dotsc,l-1$. Let $t\geq 2$. Then for any functions
\begin{equation*}
f_i\in C^{t+l-1-i}_{\mathbb{F}}(\overline{\pi_{D_i}(U)},\mathbb{F}), i=0,\dotsc,l-1,
\end{equation*}
there exists a unique $v\in C^t(\overline{U})$ such that 
\begin{equation}\label{eqn:system}
\begin{aligned}
&\mathrm{D}^2v(\cdot)(\tens{M}^{\intercal}x,y)=\mathrm{D}^2v(\cdot)(x,\tens{M}^{\intercal}y),\\
&v|_{\overline{\rho_{D_{l-1}}\pi_{D_{l-1}}(U)}}=\mu_{\mathbb{F}}f_{l-1}\circ \pi_{D_{l-1}},\\
&\mathrm{D}v(\cdot)((\tens{M}^{\intercal}-\lambda I)^{l-1-i}x)|_{\overline{\rho_{D_i}\pi_{D_i}(U)}}=\mu_{\mathbb{F}}\mathrm{D}(f_i\circ \pi_{D_i})(\cdot)(x),\\& x\in \rho_{D_i}\pi_{D_i}(B),i=0,\dotsc,l-2.
\end{aligned} 
\end{equation} 
The unique solution is given by $v=\phi_{\mathbb{F}}(\mathrm{T}f)$, where
\begin{equation*}
\mathrm{T}f=\sum_{k=0}^{l-1}\left(\rho_{(e^k)}f_k(\pi_{D_k}(u))+\sum_{j=1}^{l-1-k}\frac1{j!}\mathrm{G}^j_{(e^k)}\mathrm{D}^jf_k(\pi_{D_k}(u))((u-\rho_{D_k}\pi_{D_k}(u))^j)\right).
\end{equation*}
\end{theorem}
\begin{proof}
Let us show that $v=\phi_{\mathbb{F}}(\mathrm{T}f)$ solves the system. Since $\mathrm{T}f$ is $A$-differentiable, by Lemma \ref{L:components} it follows that for any $x,y\in \mathbb{R}^n$
\begin{equation*}
\mathrm{D}^2v(b)(\tens{M}^{\intercal}x,y)=\mathrm{D}^2v(b)(x,\tens{M}^{\intercal}y).
\end{equation*}
We shall show that the boundary conditions are satisfied. For this we use Lemma \ref{L:froboundary} and infer that
\begin{equation*}
\phi_{\mathbb{F}}(e^{l-1-i}\mathrm{T}f)|_{\overline{\rho_{D_i}\pi_{D_i}(U)}}=\mu_{\mathbb{F}}f_i\circ\pi_{D_i}
\end{equation*} 
and in particular
\begin{equation*}
v|_{\overline{\rho_{D_{l-1}}\pi_{l-1}(U)}}=\mu_{\mathbb{F}}f_{l-1}\circ\pi_{D_{l-1}}.
\end{equation*}
By $A$-linearity of derivative of $\mathrm{T}f$ we have
\begin{equation*}
\mathrm{D}v(b)(e^{l-1-i}x)=\phi_{\mathbb{F}}(\mathrm{D}\mathrm{T}f(b)(e^{l-1-i}x))=\mathrm{D}\phi_{\mathbb{F}}(e^{l-1-i}\mathrm{T}f)(b)(x)=\mu_{\mathbb{F}}\mathrm{D}(f_i\circ\pi_{D_i})(b)(x)
\end{equation*}
for all $b\in \overline{ \rho_{D_i}\pi_{D_i}(U)}$ and all $x\in \rho_{D_i}\pi_{D_i}(B)$. This means that $v$ satisfies (\ref{eqn:system}).
We shall show that the solution is unique. Assume that
\begin{equation*}
\begin{aligned}
&\mathrm{D}^2v(\cdot)(\tens{M}^{\intercal}x,y)=\mathrm{D}^2v(\cdot)(x,\tens{M}^{\intercal}y),\\
&v|_{\overline{\rho_{D_{l-1}}\pi_{l-1}(U)}}=0,\\
&\mathrm{D}v(\cdot)((\tens{M}^{\intercal}-\lambda I)^{l-1-i}x)|_{\overline{\rho_{D_i}\pi_{D_i}(U)}}=0,\\
& x\in \rho_{D_i}\pi_i(B), i=0,\dotsc,l-2.
\end{aligned} 
\end{equation*} 
We need to show that $v=0$. Since $U$ is short-path connected, from Theorem \ref{T:components} we infer that $v=\phi_{\mathbb{F}}(h)$ for some $h\in C^t_A(\overline{U},A)$. By Theorem \ref{T:structure-one}, $h=\mathrm{T}g$ for some 
\begin{equation*}
g\in\bigoplus_{k=0}^{l-1} C^{t+l-1-k}_{\mathbb{F}}(\overline{\pi_{D_k}(U)}, \mathbb{F}).
\end{equation*}
The first part of the proof shows that $\mu_{\mathbb{F}}g_{l-1}=0$ and for all $i=0,\dotsc,l-2$ functions $\mu_{\mathbb{F}}g_i$ have zero derivatives. If $\mathbb{F}=\mathbb{R}$, then, by Lemma \ref{L:froboundary}, $v=\phi_{\mathbb{R}}(\mathrm{T}g)=0$. 

Assume that $\mathbb{F}=\mathbb{C}$. Then, as $g_i$ are $\mathbb{C}$-differentiable, Cauchy--Riemann equations imply that $g_i$ have zero derivatives for $i=0,1,\dotsc,l-2$ and that $g_{l-1}$ is a real constant. Thus $\phi_{\mathbb{C}}(\mathrm{T}g_{l-1})=0$ and Lemma \ref{L:froboundary} implies that 
\begin{equation*}
\phi_{\mathbb{C}}(\mathrm{T}(g-g_{l-1}))=0.
\end{equation*}
Therefore $\phi_{\mathbb{C}}(\mathrm{T}g)=0$.
\end{proof}

\label{lastpage}

\begin{thebibliography}{10}

\bibitem{GM}
{\sc Gadea, P.M., Mu\~noz Masqu\'e, J.}
\newblock {A}-differentiability and {A}-analyticity.
\newblock {\em Proc. Amer. Math. Soc. 124}, 5 (1996), 1437--1443.

\bibitem{JAC2}
{\sc Jacobson, N.}
\newblock {\em Basic Algebra I}.
\newblock Basic Algebra. W.H. Freeman, 1985.

\bibitem{JAC}
{\sc Jacobson, N.}
\newblock {\em Basic Algebra II}.
\newblock Basic Algebra. W.H. Freeman, 1989.

\bibitem{JO}
{\sc Jodeit, M., Olver, P.J.}
\newblock On the equation {$grad f=\tens{M} grad g$}.
\newblock {\em Proc. R. Soc. Edinb. A 116}, 3-4 (1990), 341--358.

\bibitem{J1}
{\sc Jonasson, J.}
\newblock Multiplication for solutions of the equation {$grad f = \tens{M} grad
  g$}.
\newblock {\em J. Geo. Phys. 59}, 10 (2009), 1412--1430.

\bibitem{K}
{\sc Ketchum, P.W.}
\newblock Analytic functions of hypercomplex variables.
\newblock {\em Trans. Amer. Math. Soc. 30}, 4 (1928), 641--667.

\bibitem{KP}
{\sc Krantz, S.G., Parks, H.R.}
\newblock {\em A primer of real analytic functions}.
\newblock Birkh{\"a}user {A}dvanced {T}exts. Basler Lehrb{\"u}cher.
  Birkh{\"a}user Verlag, 1992.

\bibitem{L}
{\sc Lee, J.M.}
\newblock {\em Introduction to Smooth Manifolds}.
\newblock Graduate Texts in Mathematics. Springer, 2003.

\bibitem{MO}
{\sc Mazur, S., Orlicz, W.}
\newblock Grundlegende {E}igenschaften der polynomischen {O}perationen. {E}rste
  {M}itteilung.
\newblock {\em Stud. Math. 5}, 1 (1934), 50--68.

\bibitem{R}
{\sc Roscule\c{t}, M.N.}
\newblock {\em Func\c{t}ii monogene pe algebre comutative}.
\newblock Editura Academiei Republicii Socialiste Rom\^{a}nia, 1975.

\bibitem{S}
{\sc Skowro{\'n}ski, A., Yamagata, K.}
\newblock {\em Frobenius Algebras}.
\newblock No.~t. 1 in EMS textbooks in mathematics. European Mathematical
  Society, 2011.

\bibitem{W2}
{\sc Waterhouse, W.C.}
\newblock Analyzing some generalized analytic functions.
\newblock {\em Expo. Math. 10\/} (1992), 183--192.

\bibitem{W1}
{\sc Waterhouse, W.C.}
\newblock Differentiable functions on algebras and the equation {$grad (w) =
  \tens{M} grad (v)$}.
\newblock {\em Proc. R. Soc. Edinb. A 122}, 3-4 (1992), 353--361.

\bibitem{Z}
{\sc {\.Z}elazko, W.}
\newblock {\em Banach Algebras}.
\newblock Modern Analytic and Computational Methods in Science and Mathematics.
  Elsevier Pub. Co., 1973.

\end{thebibliography}
\end{document}